\documentclass{amsart}
\usepackage{graphicx,lscape,mathtools}
\usepackage{amsthm}
\usepackage{amsmath,amsfonts, amssymb}
\usepackage{MnSymbol}
\usepackage{ifthen}
\usepackage{caption}
\usepackage{subcaption}
\usepackage[justification=centering]{caption}
\usepackage[colorlinks=true]{hyperref}
\usepackage[export]{adjustbox}
\usepackage{xcolor}
\usepackage{soul}
\usepackage{quiver} 
\usepackage{tikz-cd} %commutative diagram

\newtheorem{theorem}{Theorem}[section]
\newtheorem{corollary}[theorem]{Corollary}

\newtheorem{lemma}[theorem]{Lemma}
\newtheorem{proposition}[theorem]{Proposition}

\theoremstyle{definition}
\newtheorem{definition}[theorem]{Definition}

\newtheorem{example}{Example}

\newcommand{\CC}{\mathbb{C}}
\newcommand{\DD}{\mathbb{D}}

\newcommand{\Jplus}{J_+}
\newcommand{\Jstar}{J_*}

\newcommand{\Jminus}{J_-}
\newcommand{\Jnot}{J_0}
\newcommand{\Kplus}{K_+}

\newcommand{\Kminus}{K_-}
\newcommand{\Knot}{K_0}
\newcommand{\FC}{\mathcal{U}}

\newcommand{\Rnab}{F_{n,a,b}}
\newcommand{\mapR}{F}

\newcommand{\lam}{\mathcal{L}}

\newcommand{\FClabel}[4]{\begin{pmatrix}
    #2 \curvearrowleft #1\\
    #3 \rcurvearrowright #4 \end{pmatrix}}

\DeclareMathOperator\Arg{Arg}

\title[Paper Fortune Teller changes in Julia Sets]{Paper Fortune Tellers in the combinatorial dynamics of some generalized McMullen maps with both critical orbits bounded}
\author{Suzanne Boyd and Kelsey Brouwer}
\date{\today}

\begin{document}

\begin{abstract}
    For the family of complex rational functions known as ``Generalized McMullen maps'',  
    \noindent \mbox{$\Rnab(z) = z^n + \dfrac{a}{z^n}+b$,} 
    for $a\neq 0$ and $n \geq 3$ fixed, we reveal, and provide a combinatorial model for, some new dynamical behavior. In particular, we describe a large class of maps whose Julia sets contain both infinitely many homeomorphic copies of quadratic Julia sets and infinitely many subsets  homeomorphic to a set which is obtained by starting with a quadratic Julia set, then changing a finite number of pairs of external ray landing point identifications, following an algorithm we will describe. 
\end{abstract}

\maketitle

%===========================
\section{Introduction}
\label{sec:intro}
%===========================

As a simple starting example, consider the family of quadratic polynomials $P_c(z) = z^2 + c,~c \in \mathbb{C}$. We define the \textit{Fatou set} of $P_c$ in the typical way, as the set of values in the domain where the iterates of $P_c$ is a normal family in the sense of Montel. The \textit{Julia set}, $J$, is also defined the usual way as the complement to the Fatou set. The \textit{filled Julia set}, $K$, is the union of the Julia set and the bounded Fatou components.

Douady and Hubbard (\cite{douady_hubbard}) showed that quadratic Julia sets can result from other iterative processes as well. They defined what it means for a map to behave like $P_c$, calling such a map {\em polynomial-like of degree two} (see Section \ref{sec:background}). 

Our family of interest is: $$\Rnab(z) = z^n + \dfrac{a}{z^n}+b~,~n \in \mathbb{N},~a \in \mathbb{C}^*,~b \in \mathbb{C},$$ where $\mathbb{C}^* = \mathbb{C}\backslash \lbrace 0\rbrace$. We consider $n$ a fixed integer with $n \geq 3$, thus have two complex parameters $a\neq 0$ and $b$. This family has $2n$ critical points, but only two critical values which we'll call $v_+$ and $v_-$. 

This family has been studied previously by Devaney and colleagues as well as the first author and colleagues. For $n\geq 2$, Devaney and colleagues study the subfamily with $b=0$, ``McMullen maps", in papers such as \cite{dhalo,devsurvey}. When $b=0$ the critical value orbits behave symmetrically so there is only one free critical orbit. In \cite{DevaneyRussell,devkoz}, Devaney and coauthors study the family in the case where the critical value is a fixed point. 

In \cite{devgar}, Devaney and Garijo  study Julia sets as the parameter $a$ tends to $0$, for a different generalization of McMullen Maps: $z \mapsto z^n + \dfrac{a}{z^d}$ with $n \neq d$. In \cite{dhalo,JSM2017} the authors study parameter space in the $n\neq d$ (but still $b=0$) case.

In \cite{boydschul, boyd_mitchell}, the first author and colleagues study the parameter and dynamical planes in the case $b\neq 0$, including establishing some parameters for which quadratic polynomial-like behavior exists (both in dynamical as well as parameter space). 

It is already known that this family, even for $b=0$, exhibits some behavior that is distinct from polynomial dynamics. For example, there are Julia sets which are Cantor sets of simple closed curves (\cite{mcmullen_example, devlook}, in the case that both critical values lie in the same preimage component of the basin of attraction of infinity. 
Xiao, Qiu, and Yin (\cite{XQY2014}) 
establish a topological description of the Julia sets (and Fatou components) of $\Rnab$ according to the dynamical behavior of the orbits of its free critical points. This work includes a result that if there is a critical component of the filled Julia set which is periodic while the other critical orbit escapes, then the Julia set consists of infinitely many homeomorphic copies of a quadratic Julia set, and uncountably many points. In the case of both critical orbits bounded, they show there are no Herman rings, and if both critical values are in different Fatou components, then every Fatou component is a topological disk, or if both critical values are in the same Fatou component it is infinitely connected.

In this article, we reveal another novel phenomena in this family by describing a class of maps in which there exist both infinitely many homeomorphic copies of quadratic Julia sets as well as infinitely many subsets homeomorphic to a set which is obtained by starting with a quadratic Julia set, then changing a finite number of external ray landing point identifications, combining some Fatou components while splitting others. We will describe this in detail, but the essential step is taking some sets of rays that are identified, breaking the identifications then gluing them into different pairings. This can be visualized somewhat as holding an origami paper ``fortune teller'' open in one direction, then closing it and opening it in the other.   

\begin{figure}[h]
    \centering
    \begin{minipage}{0.455\textwidth}
        \includegraphics[width=\textwidth]{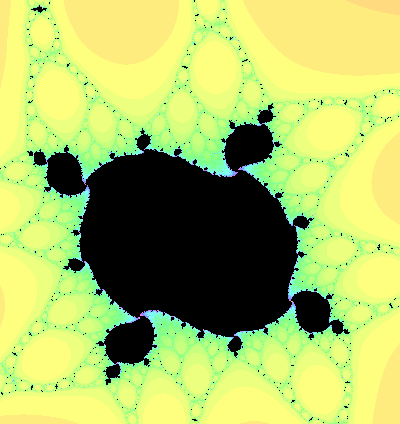}
    \end{minipage}
   % \hfill
% \hspace{3mm}
    \begin{minipage}{0.445\textwidth}
        \includegraphics[width=\textwidth]{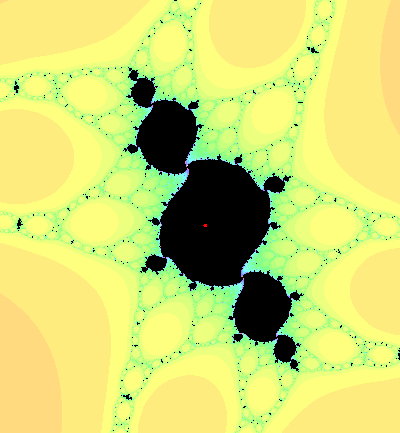}
    \end{minipage}
    \caption{\label{fig:introexample} Portions of the Julia set of $\Rnab$ for $n=5$, $a=0.1317-0.0073i$, $b=0.03+0.02i$. The right image is homeomorphic to a quadratic Julia set, it is a preimage copy under $\Rnab$ of a baby quadratic Julia set in $J(\Rnab)$ associated with the critical value $v_+$. The left is a preimage of the right under $\Rnab$, and is what we refer to as ``altered''. The red dot in the center of the right image marks the location of the other critical value $v_-$; note this is not in the location of a critical value in a basilica Julia set, which would be in one of the largest Fatou components adjacent to the central component. This alternate placement is what causes the altered shape of the preimage. }
\end{figure}
 
We call these \textbf{altered} quadratic Julia sets. See Figures~\ref{fig:introexample}, \ref{fig:exType1-2}, \ref{fig:Type2}, and~\ref{fig:Type4} for a few examples.
This occurs in the case where one critical orbit is in a homeomorphic copy of a quadratic Julia set, tending to an attracting period cycle.
The map $\mapR$ then creates an infinite tree of preimage homeomorphic copies of this homeomorphic copy. Then in our case of interest, the other critical value $v_-$ lies in one of these preimage copies; that is, it is also in the basin of attraction of the same cycle but not in the immediate basin, and not in the baby Julia set created by the dominant critical value $v_+$. However, for no reason must $v_-$ be in the ``expected'' location; that is, $v_-$ may be in any other Fatou component as compared to where the critical value would be if that copy were a quadratic Julia set. We will see that the preimage of the homeomorphic copy of the quadratic Julia set which contains $v_-$ consists of this ``altered'' set in which the ray identification changes are determined by in which Fatou set component in that preimage baby Julia set the other critical value lies. Then, this ``altered'' set has an infinite tree of preimages all homeomorphic to it.

In Section~\ref{sec:anglesdefinition} we prove our first main theorem, Theorem~\ref{thm:angle assignments}, by constructing on each baby Julia set preimage external angle assignments that respect the dynamics of $\mapR$.  In Section~\ref{sec:alteredJulias}, we provide a series of results and examples describing, using the angle assignments of Theorem~\ref{thm:angle assignments} and angle lamination diagrams, precisely how varying the position of $v_-$ within a preimage copy of a baby basilica Julia set alters the combinatorial dynamics of that preimage (and hence its preimages) so that it is no longer the same as that of a quadratic polynomial. See Propositions~\ref{prop:type0}, \ref{prop:type1M}, \ref{prop:type1notM}, Theorem~\ref{thm:typeN},
and Examples~\ref{example:Type1-1}, \ref{example:Type1-2},
\ref{example:Type2}, and~\ref{example:Type5}. We close Section~\ref{sec:alteredJulias} listing some avenues for future work.

We define external rays and angles and provide other relevant background and preliminary material, including details on our assumptions about which maps $\Rnab$ we are studying, in Section~\ref{sec:background}. 

\textbf{Acknowledgements} We thank John Hamal Hubbard for a helpful conversation about the situation of two free critical orbits, and Brian Boyd for the program Dynamics Explorer which was used to generate the images in this article.

%=========================
\section{Background and Preliminaries}
\label{sec:background}
%=========================

The generalized McMullen map $\mapR(z)=\Rnab(z)$ has $2n$ critical points, which are all possible roots of $a^{1/2n}$, but only two critical values,  $v_{\pm} = b \pm 2\sqrt{a}$. As much of the dynamics for rational functions is driven by the critical orbits, these two critical values will be vitally important to this work. We also see that $\infty$ is a superattracting fixed point, so we can define the \textbf{filled Julia set} $K= K(\Rnab)$ to be the Julia set together with all connected Fatou components whose orbits remain bounded. 

Many of the techniques utilized in this paper rely on machinery developed for polynomial dynamics, so we provide relevant background information and some assumptions in this section. In Sections~\ref{sec:anglesdefinition} and~\ref{sec:alteredJulias}, we utilize some of these techniques directly and adapt others to fit our needs.

\subsection*{External rays and angles}
For a polynomial $P$ of degree $n\geq2$ whose filled Julia set $K=K(P)$ is connected, the Böttcher coordinate $\phi: \CC \smallsetminus K \to \CC \smallsetminus \mathbb{D}$ is an isomorphism, conjugating $P$ outside of $K$ to the map $z \mapsto z^n$ outside the closed unit disk. Then the sets $\{z : \Arg(\phi(z)) \text{ is constant}\}$ generate \textbf{external rays} for $K$. Alternatively, we denote by $\rho_t$ the ray $\{\phi^{-1}(re^{2\pi it}): r>1\}$, noting that this ray is then defined on the exterior of $K$. If $\gamma(t) = \displaystyle \lim_{r\searrow1} \phi^{-1}(re^{2\pi it})$ exists, we say the ray \textbf{lands} and we can associate $z=\gamma(t)$ with the \textbf{external angle} $t$.

\begin{figure}[h]
    \includegraphics[width=0.9\textwidth]{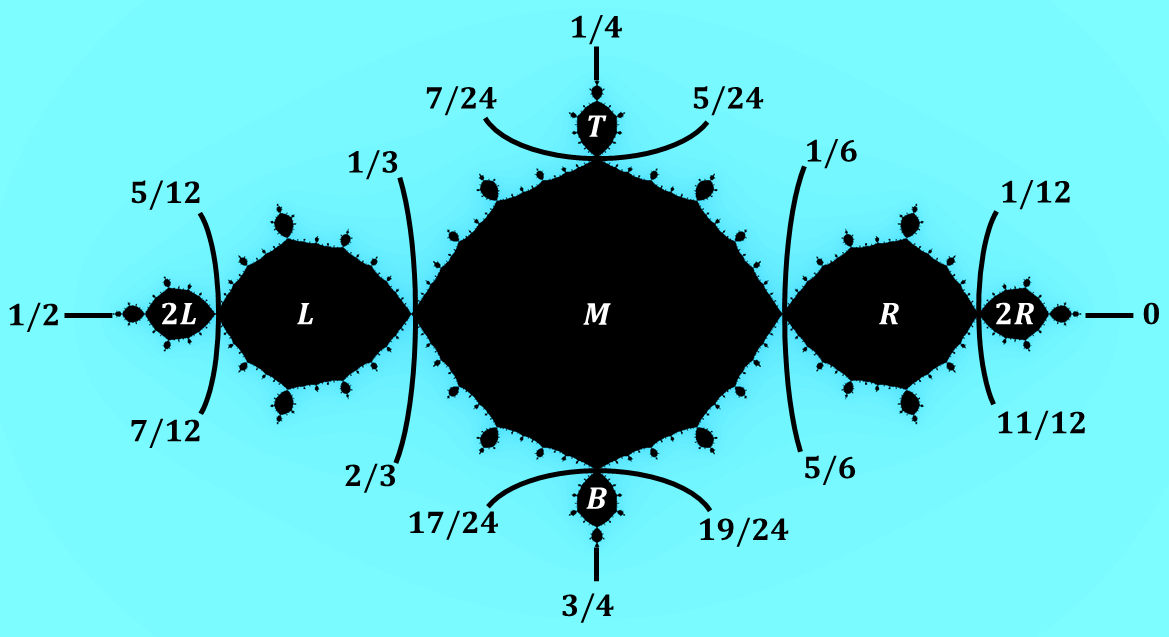}
    \caption{The Julia set for the basilica $P_{-1}(z)=z^2-1$, shown with selected external rays and letters we will use to refer to the larger Fatou components. 
    \label{fig:raysbasillicaexample}}
\end{figure}

All rays land whenever $J$ or $K$ is locally connected (see \cite{milnor}), which is true for the maps we will study. Further, if the angle $t\in \mathbb{R} \slash \mathbb{Z}$ is rational, the landing point $\gamma(t)$ is either a periodic or pre-periodic point. Hence we will focus primarily on rational angles and their landing points on a Julia set. Note that here we are using the identification of $S^1$ with the unit interval $[0,1)$, so our angles will take values in $[0,1)$.
See Figure~\ref{fig:raysbasillicaexample} for a quadratic polynomial Julia set with some labeled rays. Specifically, this is the famous ``basilica'' Julia set of the map $P_{-1}(z)=z^2-1$ for which the critical orbit lies in a (super)attracting period two cycle.

\begin{figure}[h]
    \includegraphics[width=0.9\textwidth]{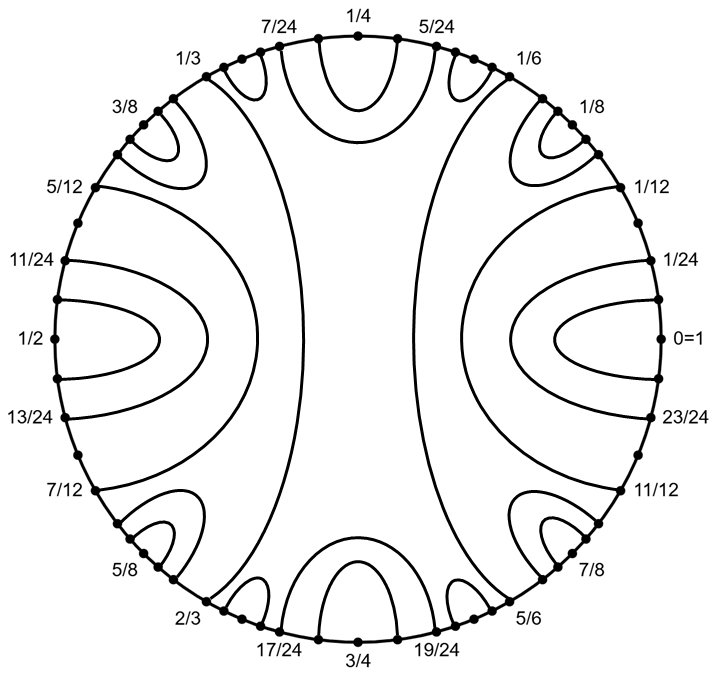}
    \caption{Lamination diagram for the basilica $P_{-1}(z)$
\label{fig:laminationexample}}
\end{figure}

Notice there are many situations in which multiple rays land at the same point. A concise way to express the external ray structure and these identifications for a Julia set is through Thurston's \cite{Thurston} \textbf{lamination diagram}. Using the $[0,1)$ parameterization of $S^1$, we connect via a \textbf{chord} any two angles $a$ and $b$ such that the rays $\rho_a$ and $\rho_b$ share a landing point, that is, if $\gamma(a)=\gamma(b)$. In this case, we will write $a \sim b$ and consider $a$ and $b$ to be \textbf{identified}. By drawing these chords so that none intersect, $\overline{\mathbb{D}}$ is partitioned into leaves which represent the components of the filled Julia set. Notably, in the case of hyperbolic rational dynamics, angles will be identified in sets of cardinality corresponding to the period of the attracting periodic orbit, or will be left unidentified. Also note for quadratic polynomials, the lamination remains invariant under the angle-doubling map. We will use $\lam$ or $\lam_{\gamma}$ to refer to the lamination of $S^1$ induced by a map $\gamma$ of $S^1$.
See Figure~\ref{fig:laminationexample}.

%-----------------------------------
\subsection*{Polynomial-like mappings}
Douady and Hubbard defined  polynomial-like maps to explain the existence of homeomorphic copies of polynomial Julia sets inside of other Julia sets. 

\begin{definition}(\cite{douady_hubbard})
\label{defn:polynlike}
A map $f:U' \to f(U')=U$ is \textbf{polynomial-like} if 
\begin{itemize}
    \item $U'$ and $U$ are bounded, open, simply connected subsets of $\CC$,
    \item $U'$ is relatively compact in $U$, and
    \item $f$ is analytic and proper.
\end{itemize}
Further, $f$ is polynomial-like of \textbf{degree two} if $f$ is a 2-to-1 map except at finitely many points, and $U'$ contains a unique critical point of $F$. The \textbf{filled Julia set of a polynomial-like map} is the set of points whose orbits remain in $U'$. 
\end{definition}

\begin{theorem} \label{thm:polynlike}
    \cite{douady_hubbard} A polynomial-like map of degree two is topologically conjugate on its filled Julia set to a quadratic polynomial on that polynomial's filled Julia set.
\end{theorem}

For this reason, we often refer to the filled Julia set of a polynomial-like map as a \textbf{baby Julia set}. 
Our results require that $\Rnab$ contains a baby Julia set.

\subsection*{Symmetry of the Julia set for generalized McMullen Maps}

The following lemma was in \cite{XQY2014}, but as the straightforward proof was not provided, we include it for benefit of the reader. 

\begin{lemma}
    For any $n \in \mathbb{N}$, $a \in \CC$, and $b \in \mathbb{C}$, the Julia set $J(\Rnab)$ and filled Julia set $K(\Rnab)$ have an $n$-fold circular symmetry.
\end{lemma}

\begin{proof}
    Let $\displaystyle x_k=re^{i(\theta+\frac{2\pi k}{n})}$, where $k=0, \dots, n-1$, be the $n$-fold circularly symmetric points to $x=re^{i\theta}$. Observe that 
    $$
\displaystyle \mapR(re^{i(\theta+\frac{2\pi k}{n}})) = r^ne^{i(n\theta + 2\pi k)} + b + \frac{a}{r^ne^{in\theta + 2\pi k)}} = r^ne^{in\theta} + b + \frac{a}{r^ne^{in\theta}} = \mapR(re^{i\theta})
$$ 
    holds for all $k$. Thus $x=re^{i\theta}$ shares an orbit with its symmetric points, and therefore all of the symmetric points must lie either in the filled Julia set $K$ or in $\CC-K$. Since $J=\partial K$, $J$ must have the same $n$-fold circular symmetry.
\end{proof}

The general symmetry gives us the following symmetry for baby Julia sets.

\begin{corollary}
\label{cor:n_baby_Js}
    If there exists one baby quadratic Julia set 
within the filled Julia set of $\Rnab$, then there are $n$ copies, with $n-1$ direct preimage copies, of the baby Julia set. Moreover, each of these $n$ sets contains a critical point of $\mapR$.
\end{corollary}

\begin{proof}
    If there exists a homeomorphic copy $K'$ of the filled Julia set $K_c$ of a quadratic polynomial 
within the filled Julia set of $\Rnab$, then the same shape must occur $n-1$ other times by the $n$-fold symmetry established in the previous lemma. Recall by Douady and Hubbard's polynomial-like map criteria, $K'$ contains a critical point of $\Rnab$. Note that the critical points of $\Rnab$ are of the form $a^{1/2n} e^{i 2\pi k/2n}$, so the $n$-fold symmetry established above places a critical point in all of the preimage copies of $K'$ as well.
    Since a baby quadratic filled Julia set can only contain one critical point, the preimage copies must be disjoint. Therefore there are a total of $n$ homeomorphic copies of $K_c$ within the Julia set of $\Rnab$.
\end{proof}

\begin{figure}[h]
        \includegraphics[width=0.6\textwidth]{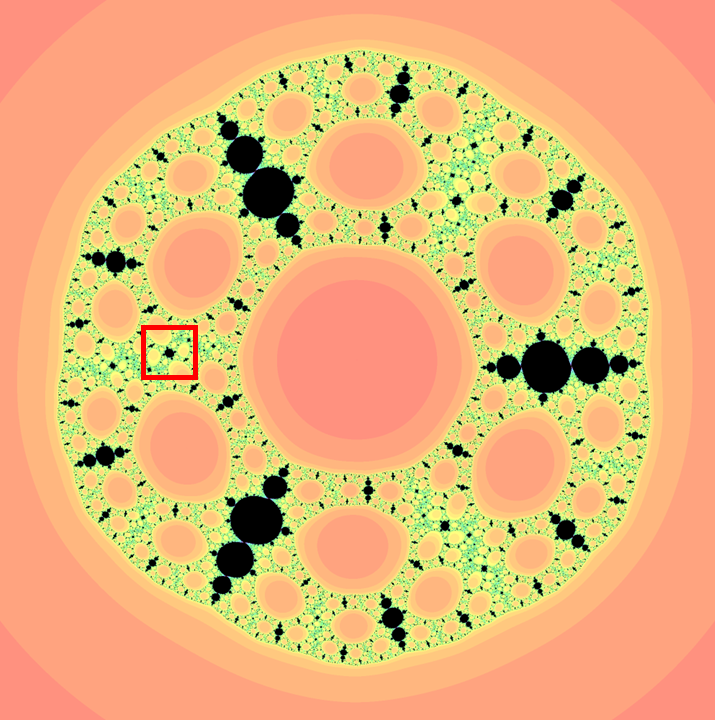}
    \caption{\label{fig:wholeJexample}
    The Julia set of $\Rnab$ for $n=3$, $a=0.05855-0.01282i$, $b=0.02+0.03i$. The baby quadratic Julia set is on the positive real axis. Its $n-1=2$ rotationally symmetric preimages are apparent, as are several smaller deeper-level preimages. The altered preimage is several preimages deep, so it is too small to see in detail without zooming in to the area outlined by the red square.}
\end{figure}

In Figure~\ref{fig:wholeJexample}, an entire filled Julia set $K(\Rnab)$ is shown. Since the altered preimage copies are several preimages deep, they are too small to see without magnification in precise places. Then, on initial inspection, one might think this Julia set was simply homeomorphic copies of a quadratic Julia set, and possibly points as in the \cite{XQY2014} case of one critical orbit in a baby $J$ and one escaping.  

\subsection*{Assumptions}
In this article, though we observed similar phenomena in many cases and are setting up machinery we expect to generalize, our main case of interest is $\Rnab$ with $n \in \mathbb{N}$, $a \in \CC^{*}$, $b \in \mathbb{C}$ and for which:
\begin{enumerate}
    \item[(A1)] $\mapR$ is hyperbolic (i.e., the parameters are in a hyperbolic component of parameter space),
    \item[(A2)] Both critical values have bounded orbits (i.e., $v_+, v_-$ of $\Rnab$ lie in the filled Julia set $K_{n,a,b} = K(\Rnab)$), but are not in the same Fatou component,
    \item[(A3)] $\mapR$ is polynomial-like of degree 2 on a region $U'$ containing one critical value (wolog $v_+$), and conjugate on the filled Julia set of $\mapR|_{U'}$, which we'll call $\Kplus$, to a quadratic polynomial $P_c$,
    \item[(A4)] $P_c$ lies in the ``basilica'' bulb; that is, $v_+$ lies in the immediate basin of an attracting cycle of period two, and finally
    
    \item[(A5)] The other critical value ($v_-$) lies in the interior of a preimage copy of the baby quadratic Julia set ($K_+$),  and therefore lies in the same attracting basin as $v_+$ but not in the immediate basin. We call that preimage $K_-$. 
\end{enumerate}

%==========================
\section{Assigning angles to copies of quadratic Julia sets in $J(\Rnab)$}
\label{sec:anglesdefinition}
%============================

Given our Assumption A3 that $\Rnab$ has a baby Julia set $\Kplus$ homeomorphic to a quadratic Julia set $K_c$, we now use the polynomial-like map to take the external angles from  $J_c = \partial K_c$ and transfer them to   $\Jplus = \partial \Kplus$, which is a subset of the filled Julia set $K(\Rnab)$. Using this and the dynamics of $\mapR$ as a guide, in this section, we will assign angles to all preimages of $\Kplus$ through backward iteration. 

These external angle assignments help us analyze the novel behavior we will present in the next section. Specifically, in Section~\ref{sec:alteredJulias}, we apply the remainder of our assumptions A1--A5 and study the case in which the second critical value is in a preimage copy of $\Kplus$; in particular, the case where the critical value contained in $P_c$ is attracted to a cycle of period 2, so that the baby Julia set $\Kplus$ takes the form of a ``basilica''. There we use the angle assignments to describe how the dynamics of the Julia sets of the rational map depend upon the location within the preimage of $\Kplus$ of the other critical value $v_-$. Generalizing these results to other quadratic baby Julia sets would be interesting future work.

In this section, we will only use Assumption A3, that $\Kplus$ is a baby quadratic Julia set. Additionally, as we want to define external angles on the preimages of $\Kplus$, we need $\Jplus$ to be connected and locally connected. Assumption A1 would imply this, but it is stronger than is needed for the following construction. 

%------------------------------------------
\textbf{Notation.}
By Assumption A3, $\mapR=\Rnab$ is polynomial-like of degree 2 on a region $U'$ containing one critical value ($v_+$), making $\mapR$ conjugate on $\Kplus$, the filled Julia set of $\mapR|_{U'}$,  to a quadratic polynomial $P_c$ on its Julia set $J_c$. Using $\Jplus = \partial \Kplus$, we denote by $\Jstar$ the infinite tree of preimages and eventual preimages of $\Jplus$; that is, let 
$$
\Jstar := \displaystyle \bigcup_{m=0}^\infty \mapR^{-m}(\Jplus).
$$ 
Next, we show how to construct external angle assignments for each point in $\Jstar$ in such a way that the local dynamics are respected.

\begin{theorem}
\label{thm:angle assignments}
    Let $\mapR=\Rnab$ be a generalized McMullen map satisfying assumption (A3), that $\mapR$ is polynomial-like on a region $U'$ containing one critical value, and is conjugate on $\Kplus$, the filled Julia set of $\mapR|_{U'}$, to a quadratic polynomial $P_c$. Further, assume the Julia set $J_c$ of $P_c$ is connected and locally connected.

    Then there exists a surjective relation $\Gamma \colon S^1 \to \Jstar$ which assigns an angle in $[0,1)=S^1$ to each point in each component of $\Jstar = \bigcup_{m=0}^\infty \mapR^{-m}(\Jplus)$ 
    so that the angle assignments respect and reflect the dynamics to and from that point. That is, if $J_{m,j} \in \mapR^{-m}(\Jplus)$ is a preimage copy of $\Kplus$, then $\Gamma$ restricted to co-domain $J_{m,j}$ is a surjective function, $\Gamma|^{J_{m,j}} \colon S^1 \to J_{m,j}$,  and: 

    \begin{enumerate}
        \item In the $m=0$ case,
        $ \mapR(\Gamma|^{\Jplus}(t)) = \Gamma|^{\Jplus}(2t).$
           
        \item If $m \geq 1$ and $K_{m,j} \neq \Kplus$ is a component of $\mapR^{-m}(\Kplus)$ which contains a critical point of $\mapR$, then
        $\mapR(\Gamma|^{J_{m,j}}(t)) = \Gamma|^{\mapR(J_{m,j})}(2t)$.

        \item 
        Finally, if $m>1$ and $K_{m,j}$ a component of $\mapR^{-m}(\Kplus)$ which does \textbf{not} contain a critical point of $\mapR$, then 
        $\mapR(\Gamma|^{J_{m,j}}(t)) = \Gamma|^{\mapR(J_{m,j})}(t)$.
    \end{enumerate}
    Practically, this means $\Gamma$ conjugates $\mapR$ to angle-doubling on the circle when the preimage copy of the baby Julia set contains a critical point; otherwise, $\Gamma$ conjugates $\mapR$ to the identity map, preserving angle assignments on the circle.     
\end{theorem}

\begin{proof}
First, we give  an overview of the proof. The construction of $\Gamma$ is inductive, beginning with assigning angles to $\Jplus$ and then assigning angles to each preimage successively. Our goal is for each eventual preimage of $\Kplus$ to have external angles assigned to each of its boundary points in a way that respects and reflects the dynamics on that preimage. In the \textbf{Step 0} below, $\Jplus$ will itself receive external angles through its identification with the quadratic Julia set $J_c$. Then we pull back angle assignments for each preimage level, where each successive pullback will fall into one of two cases. As a rational map, $\Rnab$ is conformal away from its critical points, so it is $1:1$ on any preimage that does not contain a critical point. Where $\mapR$ is $1:1$, any boundary point should map onto another point associated with the same external angles. On the other hand, a preimage may contain the other critical value, so that \textit{its} preimages contain critical points. (This is our case of interest in the next subsection, where we'll invoke assumptions (A2)-(A4) so that both critical values of $\mapR$ will lie in some (eventual) preimage of $\Kplus$, which means that all critical points of $\mapR$ also lie in preimages of $\Kplus$.) In any case, since $\mapR$ is $2:1$ nearby its critical points, $\mapR$ should act on any preimage containing a critical point by angle doubling. Therefore, the preimages of $\Jplus$ follow two cases and we assign angles to the boundary of each preimage using different rules as to reflect the dynamics happening at that preimage depending on the presence of critical points. Finally, we collect angle assignments into a single relation $\Gamma$.

As the notation $\Gamma|^{J_{m,j}}$ may be cumbersome, we  introduce  $\gamma_{m,j}\colon S^1 \to J_{m,j}$ such that $\gamma_{m,j} = \Gamma|^{J_{m,j}}$. Then each $\gamma_{m,j}$ is a function and $\Gamma$ is the relation collecting all these functions. 
On a technical note, in defining $\gamma_{m,j}$ so that it is a function, consider the question of whether two different baby Julia set preimages intersect; that is, could $J_{m,j} \cap J_{p, q} \neq \emptyset$? McMullen's book (\cite{mcmullen_book}) contains an example of such an intersection, although it wasn't for hyperbolic maps, which is our case of primary interest. If a point is in two preimage copies, then when we define the function $\gamma_{m,j}(t)=\Gamma|^{J_{m,j}}(t)=z$, we restrict to the role the point $z$ plays in the dynamics related to $J_{m,j}$, and ignore the other intersecting baby Julia set preimage, which will be assigned an angle via $\gamma_{p,q}=\Gamma|^{J_{p,q}}$. 

 \textbf{Step 0:}
Given Assumption A3, let  $\eta:  K_c \to \Kplus$ be the homeomorphism guaranteed by Theorem~\ref{thm:polynlike}, which conjugates $\mapR|_{\Kplus}$ to $P_c|_{K_c}$, i.e., $\eta\circ P_c = \mapR\circ \eta$. 
Since Fatou and Julia sets are invariant, $\eta$ restricted provides a conjugacy $\eta|_{J_c}: J_c \to \Jplus.$ Additionally, since $P_c$ is quadratic, connected, and locally connected, using the inverse of the Böttcher coordinate $\phi$, $K_c$ has external angles from $\displaystyle \gamma_c(t)=\lim_{r\searrow1}\phi^{-1}(re^{2\pi it})$. Since $P_c$ is locally connected, all rays land on $K_c$, so each point of $K_c$ has an external angle assignment.

Note that $\gamma_c$ is not $1:1$ unless the periodic orbit is a fixed point, so we expect there to be many points at which multiple external rays land, which will then be assigned multiple angles. For example, the rays of angle $1/3$ and $2/3$ both land at the same point in the $c=-1$ basilica, so $w = \gamma_{-1}(1/3) = \gamma_{-1}(2/3)$ for $P_{-1}$ where $w=(1-\sqrt{5})/2$ is the $\alpha$-fixed point. Then for some $z\in J_{m,j}$, $\Gamma^{-1}(z)$ will be more than one angle, which are the cases in which $\Gamma$ is not $1:1$.

In the $m=0$ case, we simply assign angles to $\Jplus$ using the identification to $J_c$; particularly, we construct $\Gamma|^{\Jplus}$ to assign external angles to $\Jplus$ using the conjugacies provided by $\eta$ and $\gamma_c$. We design that 
$z \in \Jplus$ is assigned the same angle as $w=\eta^{-1}(z) \in J_c$; that is, set 
$$
\Gamma|^{\Jplus}(t) :=\eta(\gamma_c(t))=z \in \Jplus.
$$
Then tracing through the conjugacies (see Figure~\ref{fig:conjugacydiagram}), 
\begin{figure}
\[\begin{tikzcd}
	{J_+} & {J+} \\
	{J_c} & {J_c} \\
	{S^1} & {S^1}
	\arrow["\mapR", from=1-1, to=1-2]
	\arrow["\eta", from=2-1, to=1-1]
	\arrow["{P_c}", from=2-1, to=2-2]
	\arrow["\eta"', from=2-2, to=1-2]
	\arrow["{\gamma_c}", from=3-1, to=2-1]
	\arrow["{z\mapsto2z}"', from=3-1, to=3-2]
	\arrow["{\gamma_c}"', from=3-2, to=2-2]
\end{tikzcd}\]
\caption{\label{fig:conjugacydiagram} Topological conjugacy diagram for Theorem~\ref{thm:angle assignments}. On $\Jplus$, the map $\mapR$ is conjugate to a $P_c$ on its Julia set $J_c$, and hence the angle assignments $\gamma_c$ for $J_c$ can be passed to $\Jplus$. }
\end{figure}
we see that 
$\gamma_c(t)=w=\eta^{-1}(z) \in J_c$, so 
$$\mapR(\Gamma|^{\Jplus}(t)) = \mapR(z) = \mapR(\eta(w)) = \eta(P_c(w)) =  \eta(P_c(\gamma_c(t)) = \eta(\gamma_c(2t)) = \Gamma|^{\Jplus}(2t).$$
Specifically, we point out that since $\mapR|_{\Jplus}$ is conjugate to a quadratic map on its Julia set, we have the invariance $\mapR(\Jplus)=\Jplus$ and the fact that each point in $\Jplus$ maps to a point with double its starting angle. Hence $\Gamma|^{\Jplus}$ conjugates $\mapR$ to angle doubling. 

We next use $\mapR^{-1}$ to pull back the angles that are assigned to $\Jplus$ and assign them to various preimages of $\Jplus$.

\textbf{Inductive Step:} Suppose angles have been assigned on all preimage components in $\mapR^{-(m-1)}(\Jplus)$, where $m\geq1$; that is, suppose we have defined a function $\gamma_{m-1,j'} = \Gamma|^{J_{m-1,j'}}: S^1 \to J_{m-1,j'}$ for each $j'$. Now, we want to assign angles to  $J_{m,j}$ for each $j$.

\textbf{Case A:} 
Suppose that $J_{m,j}=\partial K_{m,j}$ is a distinct preimage of $J_{m-1, j'}$ such that $K_{m,j}$ contains a critical point of $\mapR$. Then $J_{m,j}$ must map $2:1$ onto its image $\mapR(J_{m,j})=J_{m-1,j'}$. Thus we construct angle assignments on $J_{m,j}$ so that $\mapR$ is conjugate to angle-doubling from $J_{m,j}$ to $J_{m-1,j'}$. 
We begin assigning angles to $J_{m,j}$ by identifying the two points
$z_0 \neq z_{1/2} \in J_{m,j}$ that map to the point $w_0$ with angle assignment 0 in $J_{m-1,j'}$; that is, we locate $z_0$ and $z_{1/2}$ such that $\mapR(z_0)=\mapR(z_{1/2})=w_0=\Gamma|^{\mapR(J_{m,j})}(0)=\Gamma|^{J_{m-1,j'}}(0)$
$=\gamma_{m-1,j'}(0)$
, where $\gamma_{m-1,j'}$ is assumed to be defined in the inductive step. Then set $\gamma_{m,j}(0)=\Gamma|^{J_{m,j}}(0):=z_0$ and $\gamma_{m,j}(1/2)=\Gamma|^{J_{m,j}}(1/2):=z_{1/2}$, where the choice of which point is $z_0$ and which is $z_{1/2}$ can be made arbitrarily. 
Due to the symmetries in quadratic polynomials, making this choice allows us to set an orientation for $J_{m,j}$.
Traveling counterclockwise in $J_{m,j}$ from $z_0$ to $z_{1/2}$ will sweep out the {``upper half''} of the $S^1$ angle assignments, and we may call that path the ``upper half'' of $J_{m,j}$.
 
For any other $t\notin \{0, 1/2\}$, there are two points
$z_{t/2} \neq z_{(t+1)/2} \in J_{m,j}$ that both map to the point in $J_{m-1,j'}$ defined by $w_t = \gamma_{m-1,j'}(t)=\Gamma|^{J_{m-1,j'}}(t) =\Gamma|^{\mapR(J_{m,j})}(t) =  \mapR(z_{t/2})=\mapR(z_{(t+1)/2})$, where $\gamma_{m-1,j'}(t)=\Gamma|^{J_{m-1,j'}}(t)$ is defined by the inductive assumption. Assign $\gamma_{m,j}(\frac{t}{2})=\Gamma|^{J_{m,j}}(\frac{t}{2}):=z_{t/2}$ and $\gamma_{m,j}(\frac{t+1}{2})=\Gamma|^{J_{m,j}}(\frac{t+1}{2}):=z_{(t+1)/2}$,  where $z_{t/2}$ is the point which is encountered while traversing $J_{m,j}$ counterclockwise from $z_0$ to $z_{1/2}$. 

Note that for $z_{t/2}\in J_{m,j}$ on the upper half, we assign the same angle as the point in the same relative position on $J_{m-1,j'}$, i.e., the point $w_{t/2}$ 
which lies in $J_{m-1,j'}$ and has angle in $(0, 1/2)$. For $z_{(t+1)/2}$ in the lower half of $J_{m,j}$, similarly use the point $w_{(t+1)/2}$ 
which lies in $J_{m-1,j'}$ and has angle in $(1/2, 1)$. See Figure~\ref{fig: angle assignments j'}.

   \begin{figure}[h]
        \includegraphics[width=0.9\textwidth]{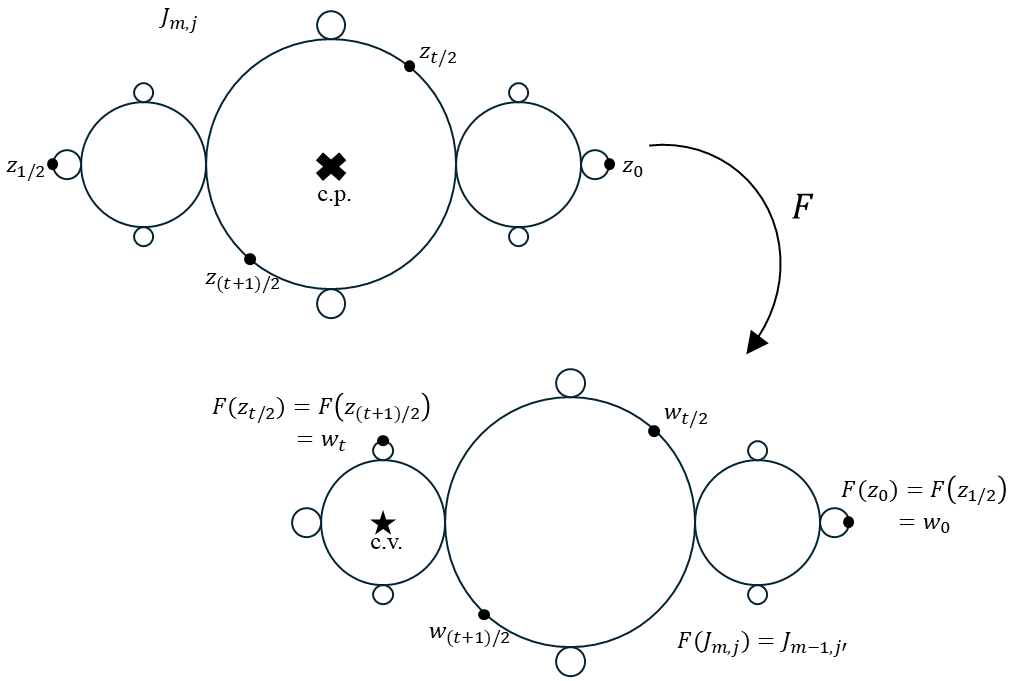}
        \caption{A diagram demonstrating angle assignments for $J_{m,j}$ if $K_{m,j}$ contains a critical point (\textbf{Case A}).
        \label{fig: angle assignments j'}}
    \end{figure}

This assigns angles to every point in $J_{m,j}$, and produces the conjugacy statements
$$
\Gamma|^{\mapR(J_{m,j})}(t) = \mapR(z_{t/2}) = \mapR(\Gamma|^{J_{m,j}}({t}/{2})) \text{ \ and } $$ 
$$\Gamma|^{\mapR(J_{m,j})}(t) = \mapR(z_{(t+1)/2}) = \mapR(\Gamma|^{J_{m,j}}({(t+1)}/{2})),$$  or equivalently, 
$\gamma_{m-1,j'}(t)=\mapR(\gamma_{m,j}(\frac{t}{2})),$ and $\gamma_{m-1,j'}(t)=\mapR(\gamma_{m,j}(\frac{t+1}{2}))$. Hence by replacing $\frac{1}{2}t$ with $t$ and $t$ with $2t$ we have
$$
\gamma_{m,j}(t) = (\mapR^{-1})|^{J_{m,j}}(\gamma_{m-1,j'}(2t))
$$
and the conjugacy from the theorem statement is established.

\textbf{Case B:} Now suppose $J_{m,j}=\partial K_{m,j}$ is a distinct preimage of $J_{m-1,j'}$ such that $K_{m,j}$ does not contain a critical point of $\mapR$. Then $\mapR|_{J_{m,j}}$ is $1:1$, 
so here each $z \in J_{m,j}$ should have the same angle as its image in $J_{m-1,j'}$, which received its angle assignments in the inductive hypothesis. Hence we define 
\begin{eqnarray*}
\gamma_{m,j}(t) =\Gamma|^{J_{m,j}}(t) & :=& 
(\mapR^{-1})|^{J_{m,j}}(\Gamma|^{\mapR(J_{m,j})}(t)) \\
& =& (\mapR^{-1})|^{J_{m,j}}(\Gamma|^{J_{m-1,j'}}(t))=(\mapR^{-1})|^{J_{m,j}}(\gamma_{m-1,j'}(t)).
\end{eqnarray*}
Then the conjugacy statement follows directly from the definition of $\gamma_{m,j}=\Gamma|^{J_{m,j}}$.

\textbf{Observation:} Recall that by assumption, $K_+$ contains one critical value of $\mapR$. Hence if $m=1$, assigning angles to $J_{m,j}=J_{1,j}$ must fall under \textbf{Case A}. If we assume that the other critical value also lies in $\Jstar$, then \textbf{Case A} occurs a second time for some $J_{\ell, j}$ where $\ell\geq1$. However these two occurrences of \textbf{Case A} account for all critical points of $\mapR$, so $J_{m,j}$ for all other values of $m (>1)$ falls into \textbf{Case B}. Thus we have defined $\Gamma$ as claimed. 
\end{proof}

Recall the angle assignments $\gamma_c$ for a polynomial map induces a lamination $\lam_c$ of the unit disk, which is invariant under angle-doubling. For the laminations $\lam_{m,j}$ induced by the $\gamma_{m,j}$'s, the previous theorem implies the following.  

\begin{corollary}    
    \label{cor:lamination}
    The lamination $\lam_{m,j}$ of the unit disk induced by any given $\gamma_{m,j}$ maps under angle-doubling to the lamination $\lam_{m-1,j'}$ induced by any $\gamma_{m-1,j'}$. 
\end{corollary}
This in turn yields:
\begin{corollary}  \label{cor:lamination_symmetry}
    For every $\gamma_{m,j}$ as defined in Theorem~\ref{thm:angle assignments},  its corresponding lamination diagram $\lam_{m,j}$ is rotationally symmetric by $180^{\circ}$.
\end{corollary}

%====================================
\section{Alterations of baby Julia sets: In which angles split and re-identify}
\label{sec:alteredJulias}
%======================================
Now that our tree of preimages has angle assignments, we can use them to describe how at a specific step in a chain of preimages of $\Jplus$, a change occurs in the angle identifications and thus in the associated lamination diagram. 
This reflects an altered shape of the baby $J$'s to bifurcate some Fatou components and combine others together, so they'll no longer look like quadratic polynomial Julia sets, and the lamination providing the combinatorial model is different from that of the quadratic polynomial Julia set. 

To clearly explain what is happening, we focus on a concrete example and suppose that $K_+$ takes the form of a baby ``basilica'', that is, that $\mapR$ is polynomial-like on $\Kplus$ and conjugate to a quadratic polynomial in the same hyperbolic component as $P_{-1}$, so $v_+$ lies in the immediate basin of an attracting period two cycle.

\subsection{Naming Fatou Components}
To understand the dynamics happening between a baby Julia set and its preimages within the Julia set of $\Rnab$, we need a convention for referring to a specific connected Fatou component lying in either the baby Julia set or one of its preimages. Per \cite{XQY2014}, in our case of interest with both critical orbits bounded but in distinct Fatou components, all connected components of the filled Julia set are topological disks. 

\textbf{Notation.}
For any connected portion of the Julia set which has been assigned external angles, we label a component $\FC$ by 
$\FClabel{a_1}{b_1}{a_2}{b_2}$
where $a_1 \sim b_2$ and $b_1 \sim a_2$ are the external angles that share landing points on the boundary of $\FC$ with lowest denominator, and where $0 \leq a_1 < b_1 < a_2 < b_2 \leq 1$. 

Consider our case of interest in which $\Kplus$ is homeomorphic to a quadratic Julia set in the basilica bulb, which is where $v_+$ is in the immediate basin of a period two cycle. We describe how to use this notation to visualize the arrangement of the angles for various Fatou components. 
For Fatou components lying along the real axis, the point with $b_1 \sim a_2$ is directly left of the point with $a_1\sim b_2$ (the latter of which is the first point of the component encountered when walking from the point assigned angle $0$ counterclockwise around the Julia set). Then we think of the component as horizontally aligned and the two identifications can be considered ``vertical pinches'' separating this Fatou component from others nearby. 

For a second case, if a Fatou component bulb is not along the main horizontal axis but in the top half, i.e., is identified by angles all of which fall in $(0,1/2)$, then tilting one's head right will line up the bulb with the matrix label.
For example, see 
$T=\FClabel{{5}/{24}}{{11}/{48}}{{13}/{48}}{{7}/{24}}$
in Figure~\ref{fig:raysbasillicaexample}.  The matrix notation provides a consistent visual if the viewer tilts their head $90^{\circ}$ to the right. 

Finally, for any bulb not along the main line but with angles in $(1/2, 1)$, one needs to tilt one's head left to get the bulb to match this matrix notation. For example, the matrix notation for $B=
\FClabel{{17}/{24}}{{35}/{48}}{{37}/{48}}{{19}/{24}}$ lines up if the viewer tilts their head $90^{\circ}$ to the left. 

We often use the label for the name of the component. For example, we  refer to the Fatou component of the filled Julia set of $P_{-1}(z)=z^2-1$ which contains the origin by $M=
\FClabel{1/6}{1/3}{2/3}{5/6}
$. This description using two pairs of identified angles holds for quadratic polynomials and can be reused for our family of functions since all local dynamics are of degree two so long as the viewer rotates their perspective to align the baby Julia set copy/preimage roughly horizontally and to the left of the point with angle assignment $0$.

We refer to two Fatou components of the filled Julia set as \textbf{adjacent} if their boundaries touch. Note that if two Fatou components share a boundary point, that point is the location of an angle identification. A simple example is that $M=\FClabel{1/6}{1/3}{2/3}{5/6}$ is adjacent to $L=
\FClabel{1/3}{5/12}{7/12}{2/3}
$ because the two components share the $\alpha$-fixed point as a boundary point, which is the point identified with both the angles ${1}/{3}$ and ${2}/{3}$. Typically, however, only one of two adjacent components will contain the identified angles of their shared boundary point in its constructed name; generally, the ``smaller'' component will. 

\begin{figure}[h]
    \includegraphics[width=0.9\textwidth]{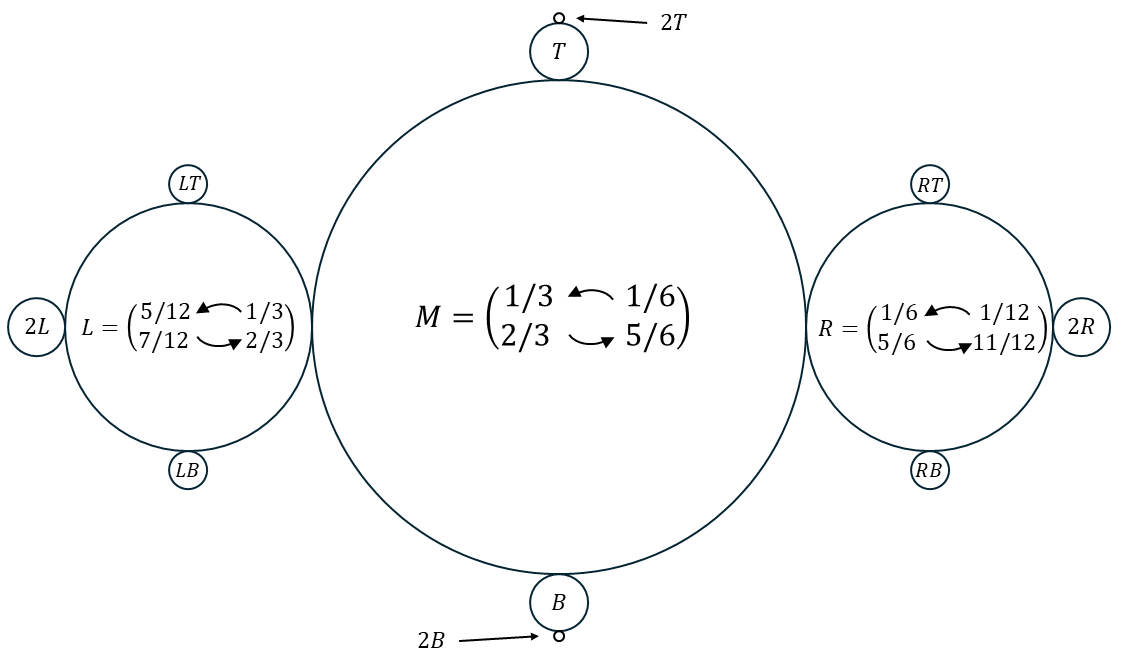}
    \caption{Diagram naming some Fatou components of the basilica Julia set}
    \label{fig:labels}
\end{figure}

Now, in the case of the baby quadratic Julia set being a basilica, referring to Figure~\ref{fig:labels}, we use the shorthand $M$ to refer to the central component, i.e., $M=\FClabel{1/6}{1/3}{2/3}{5/6}$,
and note $M$ is the location of the critical point in the basilica. We use $L$ and $R$ to refer to the largest components adjacent to $M$, so 
$L=
\FClabel{1/3}{5/12}{7/12}{2/3}
$ 
where $L$ is the location of the critical value in the basilica, and 
$R=
\FClabel{1/12}{1/6}{5/6}{11/12}
$.  
We use $2L$ for the largest component adjacent to $L$ in the direction away from $M$, and $2R$ similarly on the other side. Also, we use $T$ and $B$ for the largest components of the basilica on top and bottom, respectively, and $2T$ for the largest component on top of $T$ away from $M$ and $2B$ similarly on bottom. We also use $LT$ and $LB$ for the largest components on the top and bottom, respectively, of the component $L$, and similarly use $RT, RB$ for the largest components above and below $R$, respectively. 

\subsection{Altered baby Julia sets}
We now begin to use the full set of assumptions listed in Section~\ref{sec:background}. Notably, we now consider the assumption (A5) that while one critical value $v_+$ of $\mapR$ lies in the baby Julia set $\Kplus$ contained in $J_{n,a,b}$, the other critical value $v_-$ lies inside some preimage copy of $\Kplus$. 

\textbf{Notation:} Consistent with the notation from the prior section, suppose that $v_- \in K_{\ell, j}$, where $K_{\ell,j}$ is a distinct component of $\mapR^{-\ell}(\Kplus)$ (for some $\ell \geq 1$). 
\begin{itemize}
\item For cleaner notation, we use
$\Kminus = K_{\ell, j}$ and $\Jminus = \partial \Kminus$, 
so  $v_- \in \Kminus$ parallels $v_+ \in \Kplus$.
\item 
Let $\FC_-$ refer to the Fatou component  of $\Kminus$ in which $v_-$ lies. 
\item Let $\FC_c$ be the Fatou component of $K_c$ that uses the same identified angle name as $\FC_-$.  
\item Again for cleaner notation, we use $\Jnot=J_{\ell+1, j'}$ to refer to any preimage component in $\mapR^{-1}(\Jminus)$, so that $\mapR(J_0)= J_-$. 
\item Given a $\Jnot=J_{\ell+1,j'}$, let 
$\gamma_- = \gamma_{\ell,j}:S^1 \to \Jminus$ and 
$\gamma_0 = \gamma_{\ell+1,j'}:S^1 \to \Jnot$ be the shorthand notation for the $\Gamma|^{\Jminus}$ and $\Gamma|^{\Jnot}$ constructed in Theorem~\ref{thm:angle assignments}.
\item Recall $\lam_\gamma$ refers to a lamination of $\DD$ induced by a map $\gamma$ of angle assignments from $S^1$ to a Julia set. Similarly to the above, we let $\lam_c$ denote the lamination associated with $\gamma_c$, $\lam_{m,j}$ for $\gamma_{m,j}$, $\lam_0$ for $\gamma_0$, $\lam_-$ for $\gamma_-$, etc. 
\end{itemize}

Although $\Kminus$ is a ``clone'' of a baby Julia set $\Kplus$, as in it has the same identified angles and relative size and position of Fatou components, nothing is forcing $v_-$ to lie in the same relative position in $\Kminus$ as $v_+$ does in $K_+$ (which is the same as $c$ does in $P_c$). A priori, $v_-$ could be anywhere in $\Kminus$, as the two critical orbits are independent.   On the right side of Figure~\ref{fig:introexample} is an example where the critical point $v_-$ is shown to be in the central Fatou component of a basilica map rather than the Fatou component on the ``left'' of center.

We find that the position of $v_-$ within $\Kminus$ determines the changes to angle identifications and Fatou components in $\Jnot$ and its tree of preimages as compared to $\Jminus$ and its copies.  We describe each possible case in detail in this subsection.

\textbf{Wolog.} Note that in this section, we always assume that the critical value $v_-$ lies somewhere on the ``left'' side of $\Kminus$ or that it lies in some smaller decoration attached to $M$. We can make this assumption because if we ever find that $v_-$ lies on the ``right'' side of $\Kminus$, we can return to the step of the angle assignment construction in which $J_{1,j}$ receives its angles and reverse the choices of ``top'' and ``bottom'' so that $v_-$ lies on the left instead. This simplifies the following analysis into fewer cases.

\begin{proposition} \label{prop:type0}
    \textbf{(Type 0)} 
    If $v_-$ lies in the Fatou component of $\Kminus$ referred to by the identified angles $L = \FClabel{1/6}{1/3}{2/3}{5/6}$, then
    $\lam_0 = \lam_- = \lam_c$; that is,
    the lamination $\lam_0$ induced by $\gamma_0$ is the same as the one $\lam_-$ induced by $\gamma_-$, which is the same as the one for $\gamma_c$. 
\end{proposition}

\begin{proof}
This is the simplest case because $v_-$ lies in the ``expected'' location which is a large Fatou component immediately to the side of $M$ in $\Kminus$, which as mentioned above can be assumed to be the ``left'' side of $M$. That is,
if $\FC_-$ and $\FC_c$ are both 
$L=\FClabel{1/3}{5/12}{7/12}{2/3}$,
then the fact that the preimage of $\FC_c$ is a single component of $K_c$ which contains the critical point of $P_c$ and maps $2:1$ onto $\FC_c$ means nothing needs to change in $\Jnot$ as compared to how it is in $\Jminus$. 
Intuitively, since the basilica already expects to have a central critical point component mapping $2:1$ onto one of the side components, $\Jnot$ has the needed structure where a component containing a critical point is mapping $2:1$ onto its image, and no changes in angle identifications are needed to create this structure as they will be in future cases.

Formally, we say that if $\FC_c$ is
$\FClabel{1/3}{5/12}{7/12}{2/3}
$, we declare $v_-$ to be in its ``expected'' position and no angles split and are re-identified, so $\Jnot$ has all of the same identified angles as $\Jminus$ and hence the same lamination, which is still the basilica lamination shown in Figure~\ref{fig:laminationexample}.

Further, we see that all preimages $\mapR^{-1}(K_-)$ appear as ``typical'' quasi-conformal copies of the basilica preserving the relative sizes of different Fatou components as in the quadratic basilica.  In this case, zooms of Julia set images show only similar looking baby basilicas.
This is the only case which doesn't lead to an altered baby Julia set shape in the Julia set of $\mapR$.
\end{proof}

\smallskip

Next, if $v_-$ is in a Fatou component other than $L$, $\FC_c$ is not $L$ so it has two preimage components in $K_c$.  However, since $\FC_-$ contains a critical value of $\mapR$, it must have a single preimage component within each $\Jnot$ which contains a critical point and hence maps $2:1$ onto $\FC_-$. This is what causes an ``alteration'' in the lamination and in the shape of the filled Julia set: the dynamics on each $\Jnot$ must change from the ones on $\Jminus$ so that the two preimage components of $\FC_c$ have their corresponding components in $\Jnot$ combined into a single component, which will contain a critical point and map $2:1$ onto $\FC_-$.

Since the Fatou component $M$ contains the critical point in the basilica, if $v_-$ lies in $M$ the alteration is a special case, so we examine that first. 

\smallskip

\begin{proposition}
     \label{prop:type1M}
    \textbf{(Type 1-1)}  
    If $v_-$ lies in the Fatou component $\FC_-$ of $\Kminus$ which uses the identified angle name $M=\FClabel{1/6}{1/3}{2/3}{5/6}$, then the lamination $\lam_0$ for $\gamma_0$ only differs from $\lam_-$ induced by $\gamma_-$ (and $\gamma_c$) in that 
the  identifications $\frac{1}{6} \sim \frac{5}{6}$ and $\frac{1}{3} \sim \frac{2}{3}$ which occur in $\lam_-$ are replaced in $\lam_0$ by the identifications $\frac{1}{6} \sim \frac{1}{3}$ and $\frac{2}{3} \sim \frac{5}{6}$. 

This results in three leaf difference between $\lam_0$ and $\lam_-$: the two components corresponding to $L$ and $R$ in $\lam_-$ are combined into one in $\lam_0$, $\FClabel{1/12}{5/12}{7/12}{11/12}$, and doing so splits the central leaf (the Fatou component corresponding to $M$) in $\lam_-$ into two in $\lam_0$: $\FClabel{1/6}{5/24}{7/24}{1/3}$
    and 
    $\FClabel{2/3}{17/24}{19/24}{5/6}$. 
\end{proposition}

\begin{figure}[h]
    \centering
    \begin{minipage}{0.45\textwidth}
        \includegraphics[width=\textwidth]{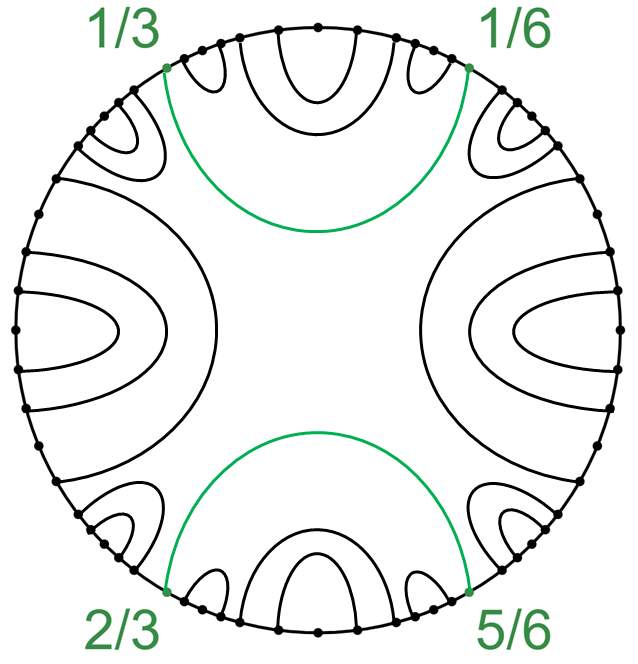}
    \end{minipage}
   % \hfill
%    \hspace{3mm}
    \begin{minipage}{0.45\textwidth}
        \includegraphics[width=\textwidth]{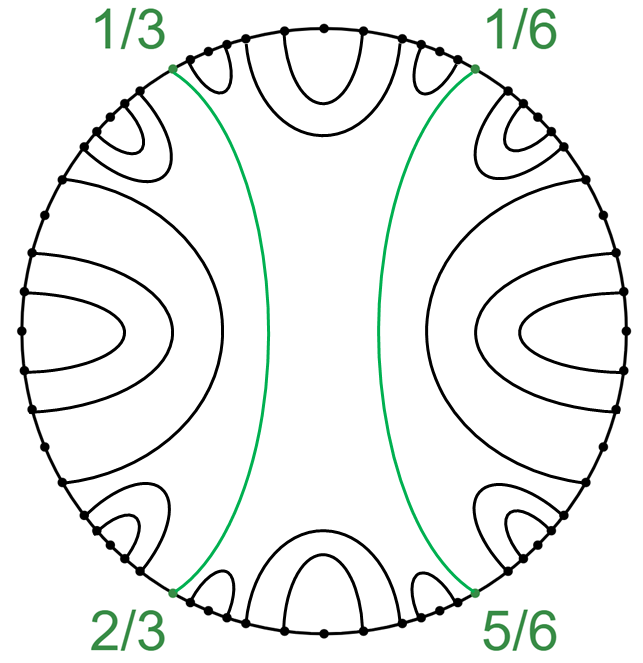}
    \end{minipage}
\caption{\label{fig:Mlamination}  \textbf{Type 1-1:} Lamination $\lam_0$ for $\gamma_0$ when $v_- \in \FC_- = M$ (left), Lamination $\lam_-$ for $\gamma_-$ is the standard basilica (right)} 
\end{figure}

\begin{figure}[h]
    \centering
    \begin{minipage}{0.4\textwidth}
        \includegraphics[width=\textwidth]{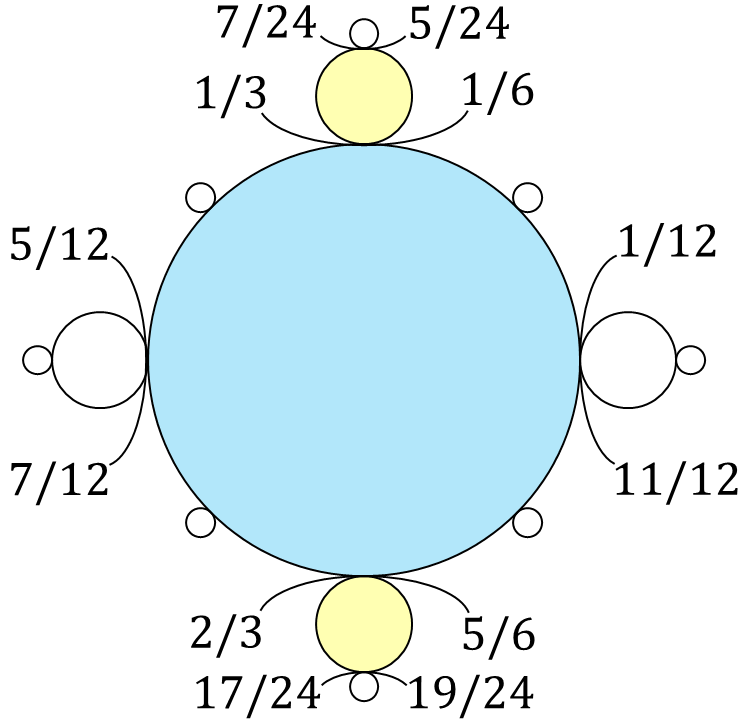}
    \end{minipage}
   % \hfill
   \hskip3mm
    \begin{minipage}{0.55\textwidth}
        \includegraphics[width=\textwidth]{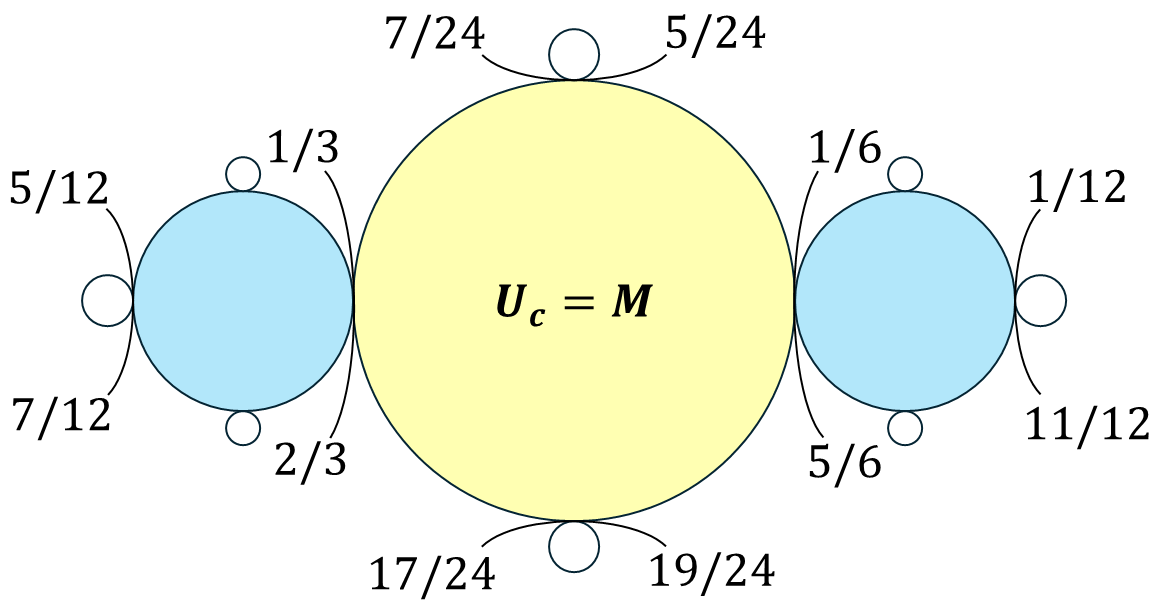}
    \end{minipage}
\caption{\label{fig:Type1-1} \textbf{Type 1-1:} When $v_- \in M$ of $\Jminus$, altered shape of $\Jnot$ (left) compared with $\Jminus$ (right). Colors have been added to show which components in $\Jminus$ are split or combined in $\Jnot$ under the new angle identifications.
} 
\end{figure}

\begin{proof} 

    Because $v_-$ is a critical value in $\Kminus$, $\Knot$ must contain a critical point.
    Then by Theorem~\ref{thm:angle assignments}, the action of $\mapR$ is conjugate to angle doubling on the angle assignments on $\Jnot$.
    Observe that $\FC_c=M$ has two preimage components in $K_c$ under $P_c$, which are $L=\FClabel{1/12}{1/6}{5/6}{11/12}$ and $R=\FClabel{1/3}{5/12}{7/12}{2/3}$. Since $\Jnot$ needs to have a single critical point component that maps $2:1$ onto $\FC_-$, the components of $\Jnot$ that correspond to $L$ and $R$ must be combined into one, and the single component $M$ of $\Jnot$ must in fact be split into two components that will both map onto $L$ to maintain dynamical structure. 

    Since the ``expected'' location of $v_-$ in $\Kminus$ is $L$ but it is actually in $M$, the point of bifurcation is the point where $L$ meets $M$, which on $\Jminus$ is the point with the identified angles $\frac{1}{3} \sim \frac{2}{3}$. 
    Since $\mapR$ has local degree two, each of $\frac{1}{3}$ and $\frac{2}{3}$ has two preimages under angle doubling; explicitly, $\frac{1}{6}$ and $\frac{2}{3}$ both double to $\frac{1}{3}$, whereas $\frac{1}{3}$ and $\frac{5}{6}$ both double to $\frac{2}{3}$. Each of these angles is identified with one from the other set on $J_c$, but on $J_c$ these identifications create two distinct components which both map onto $M$. To combine those two components, we split and re-identify these angles into the opposite pairings, which creates a single component with two points opposite each other that will map onto the same point, creating the $2:1$ effect as desired.
    Imagine holding an origami paper ``fortune teller'' open in one orientation, then closing it and opening it in the other orientation. This process is parallel to bringing the two ray pairings together at one point, and then pulling them apart with different ray landing points paired up. 
    
    In this case, we see that $\Jminus$ has $\frac{1}{6} \sim \frac{5}{6}$ and $\frac{1}{3} \sim \frac{2}{3}$ (which is the same as on $J_c$), but in $\Jnot$ the pairings must be $\frac{1}{6} \sim \frac{1}{3}$ and $\frac{2}{3} \sim \frac{5}{6}$. All other identified angles on $\Jnot$ are the same as they are on $\Jminus$, as this also respects the angle doubling established by Theorem~\ref{thm:angle assignments}. 
    See Figure~\ref{fig:Mlamination} for a side-by-side comparison of the lamination $\lam_0$ for $\gamma_0$ (left) with $\lam_-=\lam_{-1}$ (right), where the changes are highlighted in green.

    These changes result in the new critical point component in $\Jnot$ being identified by the angles 
    $\FClabel{1/12}{5/12}{7/12}{11/12}$,
    and the component which on $J_c$ maps $2:1$ onto $L$ must in $J_0$ be split into two new components, which happens naturally with the reidentifications listed above. 
    We note that this name is consistent with our naming convention: recall we label a Fatou component $\FC$ by 
$\FClabel{a_1}{b_1}{a_2}{b_2}$
where $a_1 \sim b_2$ and $b_1 \sim a_2$ are the external angles that share landing points on the boundary of $\FC$ with lowest denominator, and where $0 \leq a_1 < b_1 < a_2 < b_2 \leq 1$. Now, there are angles on the boundary of this new critical point component with lower denominators, namely $1/3, 2/3, 1/6,$ and $5/6$; however, due to their new relative positioning, these angles don't satisfy the requirements for our naming convention.

    The other new components in $\Jnot$ are 
    $\FClabel{1/6}{5/24}{7/24}{1/3}$
    and 
    $\FClabel{2/3}{17/24}{19/24}{5/6}$.
    See Figure~\ref{fig:Type1-1} for a side-by-side comparison of $\Jnot$ and $\Jminus$ showing the changes in angle identification and components. 
    
\end{proof}

\begin{example}
\label{example:Type1-1}
    Figure~\ref{fig:introexample} 
    shows computer generated images of portions of $J(\Rnab)$ with $n=5$, $a=0.0137-0.0073i$, and $b=0.03+0.02i$, with the portion around $\Jnot$ on the left and $\Jminus$ on the right. This map appears to be of \textbf{Type 1-1}, the left image maps onto the right which appears to be a typical quasi-conformal basilica. 
\end{example}

Next, while the concept explained in the proof above is consistent for all $v_-$ placements that are adjacent to $L=\FClabel{1/3}{5/12}{7/12}{2/3}$, any case where $\FC_-$ is not $M$ follows a consistent pattern in terms of which angles split and reidentify to change the lamination diagram.

\begin{proposition}
    \label{prop:type1notM}
    \textbf{(Type 1-2)} 
    Suppose that $v_-$ lies in a component $\FC_-$ of $\Kminus$ which shares an identified angle name with the component $\FC_c$ of $K_c$, and $\FC_c$ is adjacent to 
    $L=
    \FClabel{1/3}{5/12}{7/12}{2/3}$
    but is not equal to $M=\FClabel{1/6}{1/3}{2/3}{5/6}$. 
    Refer to the two preimage components of $\FC_c$ as $
    \FClabel{a_1}{b_1}{a_2}{b_2}$ and $ \FClabel{a_3}{b_3}{a_4}{b_4}$.

    Then the lamination $\lam_0$ for $\gamma_0$ differs from the $\lam_-$ of $\gamma_-$ and $\gamma_c$ only in that $\lam_-$ includes the identifications $a_1 \sim b_2$ and $a_3 \sim b_4$, which in $\lam_0$ are replaced by the identifications $a_1 \sim b_4$ and $b_2 \sim a_3$. 

    This results in three leaf differences between the laminations $\lam_0$ and $\lam_-$: two components of $\lam_-$ combine into a new component in $\lam_0$: $\FClabel{a_1}{b_2}{a_3}{b_4}$, while one component of $\lam_-$ splits in two in $\lam_-$, now $\FClabel{b_2}{1/3}{2/3}{a_3}$ and 
    $\FClabel{1/6}{a_1}{b_4}{5/6}$.

\end{proposition}

\begin{proof}
    In the same way as in Proposition~\ref{prop:type1M}, 
    for $v_-$ in a Fatou component $\FC_c$ adjacent to $L$, we have that the preimage Fatou components of $\FC_c$, labeled 
    $
    \FClabel{a_1}{b_2}{a_2}{b_2}
    $ and $
    \FClabel{a_3}{b_3}{a_4}{b_4}
    $, are adjacent to $M=
    \FClabel{1/6}{1/3}{2/3}{5/6}
    $, this time assuming that $\FC_-$ and $\FC_c$ are not $M$ itself.
    As before, by Theorem~\ref{thm:angle assignments}, the preimage of $\FC_-$ in each $\Knot$ must be one component which maps $2:1$ onto $\FC_-$ to preserve the dynamics on angle doubling nearby the critical point. Therefore the two preimage components of $\FC_c$ listed above will need to have their corresponding components in $\Knot$ combined into a single preimage component containing the critical point.
    Following the argument from Proposition~\ref{prop:type1M}, we observe that if $a<b$ are the two angles which are identified at the point where $\FC_-$ meets $L$, then $a_1$ and $a_3$ are the two angles that double to $a$ and $b_2$ and $b_4$ are the two angles that double to $b$. Recall as well that by our naming convention for the preimage components we have that $a_1 \sim b_2$ and $a_3 \sim b_4$ in $\Jminus$. As before, we can combine these two preimage components into a single component in $\Knot$ containing the critical point which maps $2:1$ onto $\FC_-$ if
    the identifications $a_1 \sim b_2$ and $a_3 \sim b_4$ are broken and reidentified as $a_1 \sim b_4$ and $b_2 \sim a_3$.
    Notably, these new identifications still reflect that these points map to the point where $\FC_-$ meets $L$ in $\Jminus$ with both identified angles at that point being represented. Additionally, this gives us that the new component in $\Knot$ containing the critical point would be named 
    $\FClabel{a_1}{b_2}{a_3}{b_4}$. 
    Meanwhile, the component in $\Knot$ that formerly was $M$ no longer contains a critical point so is split into two components which will each map onto $L$ to maintain the $2:1$ mapping of $\Knot$ onto $\Kminus$. These new components are named by the identified angles
    $\FClabel{b_2}{1/3}{2/3}{a_3}$
    and 
    $\FClabel{1/6}{a_1}{b_4}{5/6}$
 following naturally from the angle identification changes made above.
    Note again that the identifications $b_1 \sim a_2$ and $b_3 \sim a_4$ from $\Jminus$ remain identified in $\Jnot$, but these identified points now each lie on the boundary of the new central component.
\end{proof}

\medskip

\begin{figure}[h]
    \includegraphics[width=\textwidth]{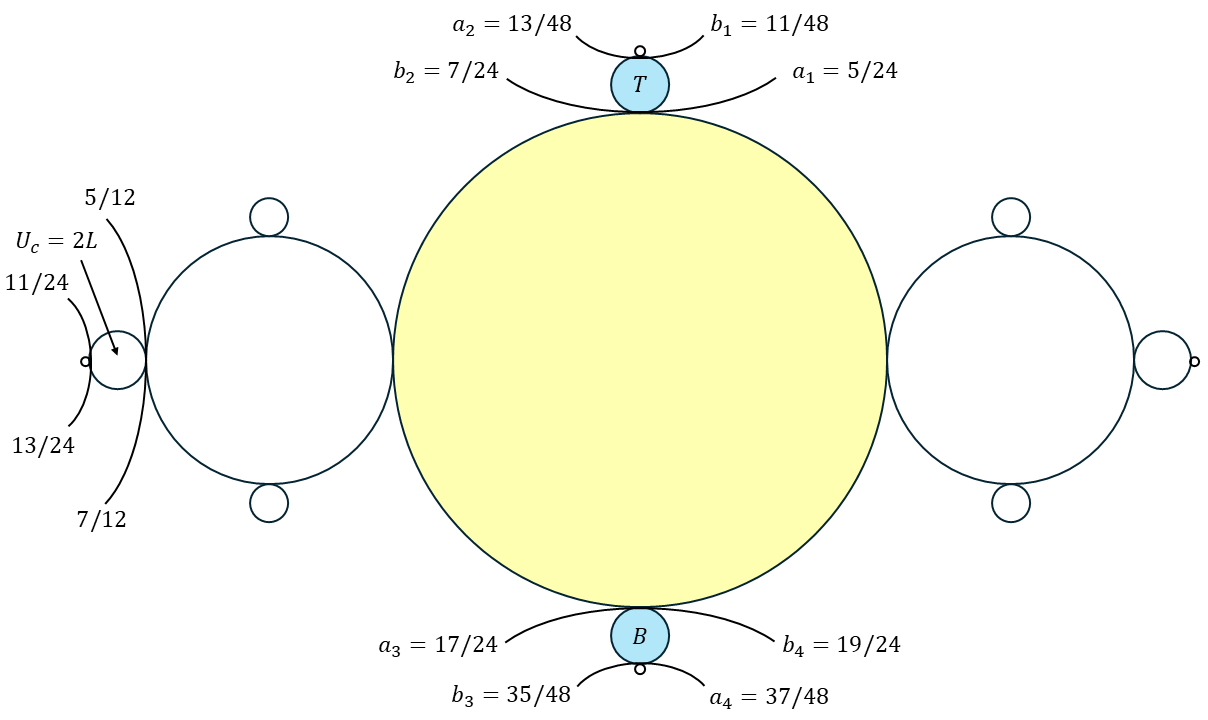}
    \caption{\textbf{Type 1-2:} Angles labeled on $J_c$ where $\FC_c = 2L$. The same angles are identified on $\Jminus$. When pulled back to $\Jnot$, the yellow Fatou component will split in two and the blue Fatou components will combine as shown in Figure~\ref{fig:2Lreid}.
    \label{fig:2Lexample}}
\end{figure}

\begin{example}
\label{example:Type1-2}
For illustrative purposes, consider the case of \textbf{Type 1-2} where $\FC_-=\FC_c= 
2L =
\FClabel{5/12}{11/24}{13/24}{7/12}
$. 
Then the preimages under $P_c$ of $\FC_c$ are 
$\FClabel{a_1}{b_1}{a_2}{b_2}
= T=
\FClabel{5/24}{11/48}{13/48}{7/24}
$ and $
\FClabel{a_3}{b_3}{a_4}{b_4}
= B =
\FClabel{17/24}{35/48}{37/48}{19/24}
$, where we see that 
$\frac{5}{24} \sim \frac{7}{24}$, $\frac{11}{48} \sim \frac{13}{48}$, $\frac{17}{24} \sim \frac{19}{24}$, and $\frac{35}{48} \sim \frac{37}{48}$.
See Figure~\ref{fig:2Lexample}.
\end{example}

We have that each identification shares a common denominator and that all of $a_1, b_2, a_3,$ and $b_4$ have the same denominator, which is smaller than the shared denominator of the other four angles, $b_1, a_2, b_3$, and $a_4$. 
Note that these are the identifications on $\Jminus$, which are the same as the typical expected basilica identifications. 
However, since $\FC_-$ contains a critical value, its preimages should each be a single component containing a critical point that maps $2:1$ onto $\FC_-$, so following the second point of Theorem~\ref{thm:angle assignments},
the angles demarking the two preimage components of $\FC_c$ must be split and reidentified in such a way that combines these two components into a single component which contains the critical point. Further, since $v_-$ placements of Type 1-2 are assumed to be in a component which is adjacent to $L$ other than $M$,
and one can observe that $M$ is the sole preimage of $L$ in $K_c$, the preimages of $\FC_c$ must lie adjacent to $M$.

\begin{figure}[h]
    \includegraphics[width=\textwidth]{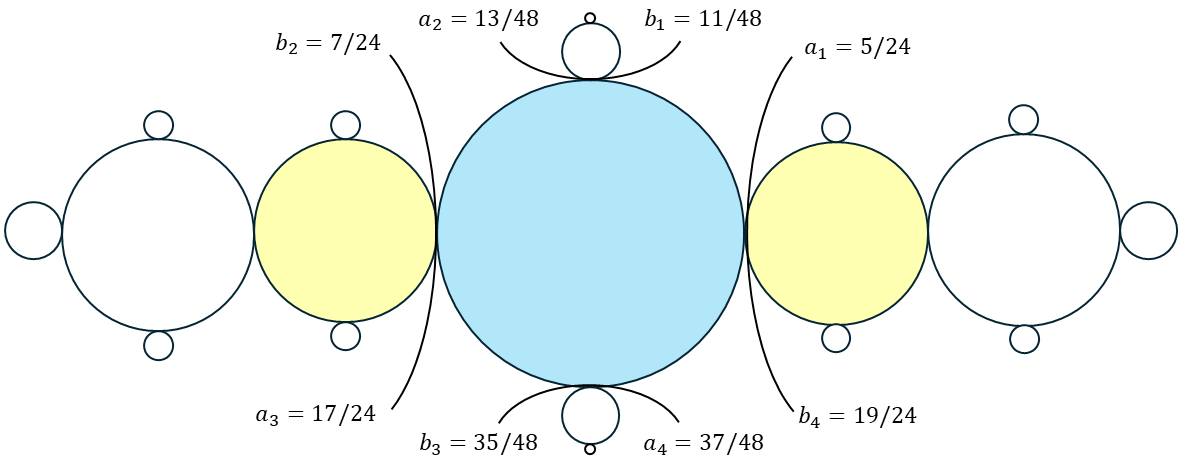}
    \caption{\textbf{Type 1-2:} Angles labeled on $\Jnot$ when $v_- \in \FC_- = 2L$. As in Figure~\ref{fig:2Lexample}, the yellow components are the result of the yellow Fatou component from $\Jminus$ splitting in two,  the blue Fatou component is the two blue components from $\Jminus$ merged into one.    \label{fig:2Lreid}}

\end{figure}

\begin{figure}[h]
    \centering
    \begin{minipage}{0.45\textwidth}
        \includegraphics[width=\textwidth]{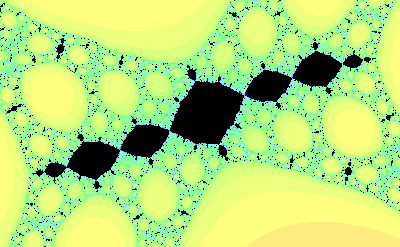}
    \end{minipage}
   % \hfill
\hspace{3mm}
    \begin{minipage}{0.45\textwidth}
        \includegraphics[width=\textwidth]{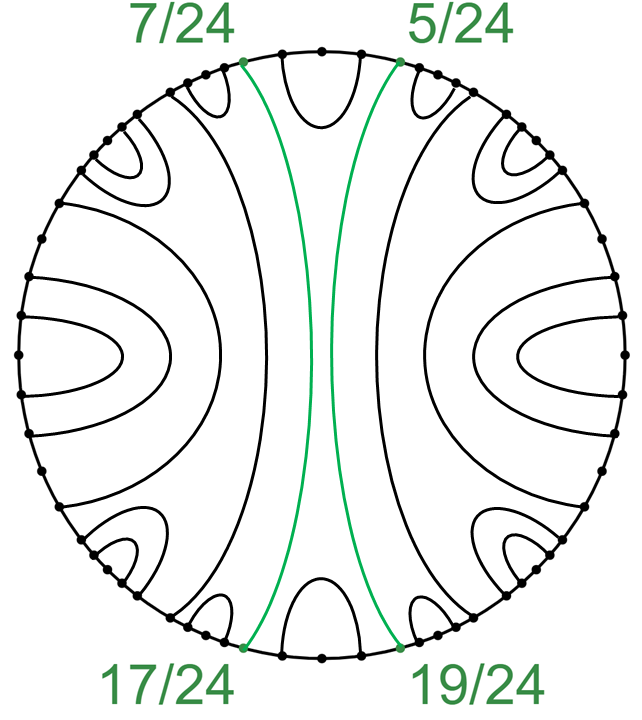}
    \end{minipage}
\caption{\label{fig:exType1-2} \textbf{Type 1-2 example:} Altered shape of $\Jnot$ when $v_- \in \FC_- = 2L$ of $\Jminus$ (left, $n=3$, $a=0.0522-0.01292i$, $b=0.01+0.03i$) along with its lamination diagram 
$\lam_0$ showing the identified angles on $\Jnot$ (right). 
} 
\end{figure}

Hence on $\Jnot$, the $a_1 \sim b_2$ identification is broken, as well as the $a_3 \sim b_4$ and instead we get these two pair identifications changed to $a_1 \sim b_4$ and $b_2 \sim a_3$, and all the other angle identifications from $\Jminus$ are unchanged including $b_1 \sim a_2$ and $b_3 \sim a_4$. So, $T$ and $B$ combine, while $M$ is split in two components on either side of the new combo component. 
See Figure~\ref{fig:2Lreid} for the key angles labeled on $\Jnot$ for this example of $v_- \in 2L$, and see Figure~\ref{fig:exType1-2} for a computer generated image of the $\Jnot$ portion of a $J(\Rnab)$, alongside the lamination of $S^1$ induced by this $\gamma_0$.

Note that Type 1-1 is only distinct from Type 1-2 in that the angles that split and re-identify do not follow the same ordering and common denominator rules. The changes in identifications that happen in both Type 1-1 and Type 1-2 cases alter the shape of $\Jnot$ so that it no longer appears as a clone of its image $\Jminus$. When the two preimage components of $\FC_-$ are combined into one, that component will now visually lie in the center of $\Jnot$, and the reidentification also squeezes the previous central component of $\Jnot$ into two components. However, in these new components there are not any critical points, so $\mapR$ acts conformally, and we see that any tree of neighboring components moves with these new split components. This allows us to understand the shape of $\Jnot$ from its angle identifications, along with the size and relative position of its components. 

For our most general result, we will explain how when $v_-$ is more than one Fatou component away from the ``expected'' location of $L$, we can describe the alterations in $\Jnot$ and $\gamma_0$ by a chain of alterations similar to those of Type 1. 

However, to help introduce the idea of performing multiple lamination changes, we first explain an example requiring only two steps. 

\begin{figure}[h]
    \centering
    \begin{minipage}{0.32\textwidth}
        \includegraphics[width=\textwidth]{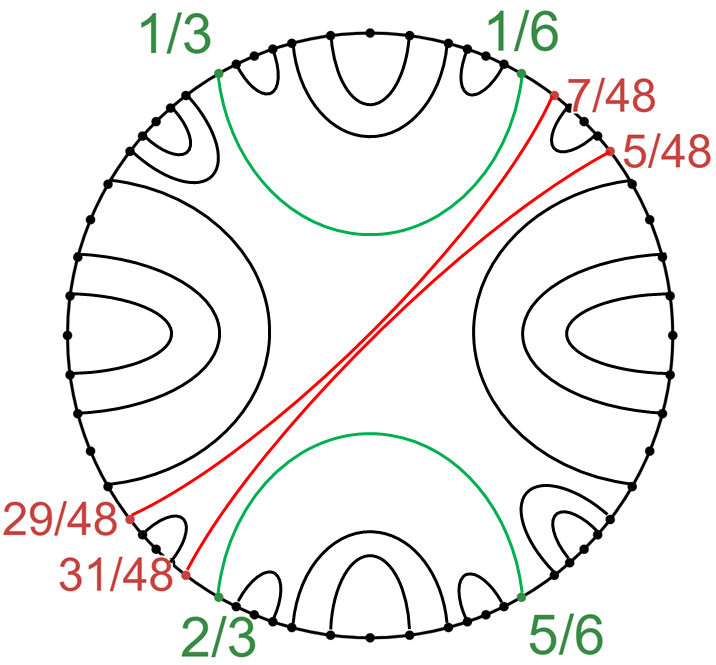}
    \end{minipage}
    \begin{minipage}{0.32\textwidth}
        \includegraphics[width=\textwidth]{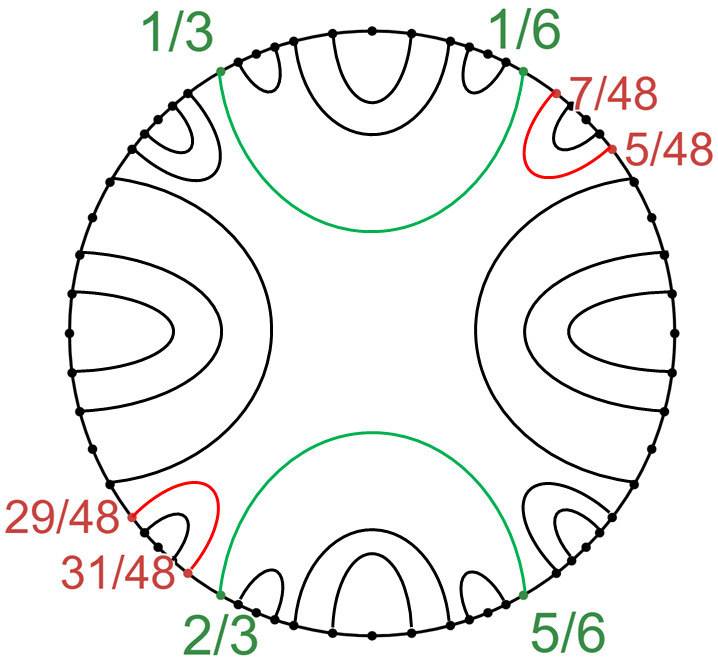}
    \end{minipage}
    \begin{minipage}{0.32\textwidth}
        \includegraphics[width=\textwidth]{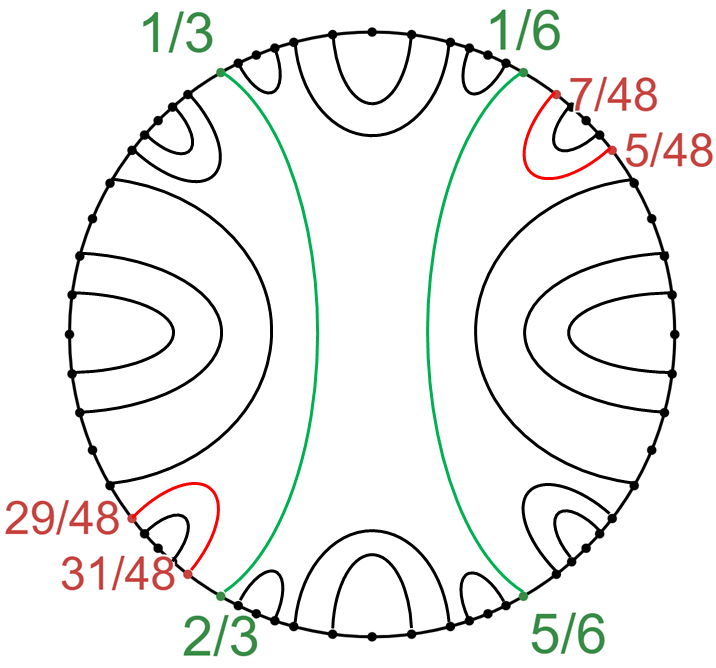}
    \end{minipage}
\caption{\label{fig:type2lamsteps} Chain of two steps of lamination changes for a \textbf{Type $N=2$} example, when $v_-\in T$. On the right is $\lam_-$ with the angles that need to break and re-identify shown in red, and the obstructing angles shown in green. In the middle is the intermediate lamination where the green angles are changed. On the left the red angles are reidentified in the lamination $\lam_0$. }
\end{figure}

\begin{figure}[h]
    \centering
    \begin{minipage}{0.3\textwidth}
        \includegraphics[width=\textwidth]{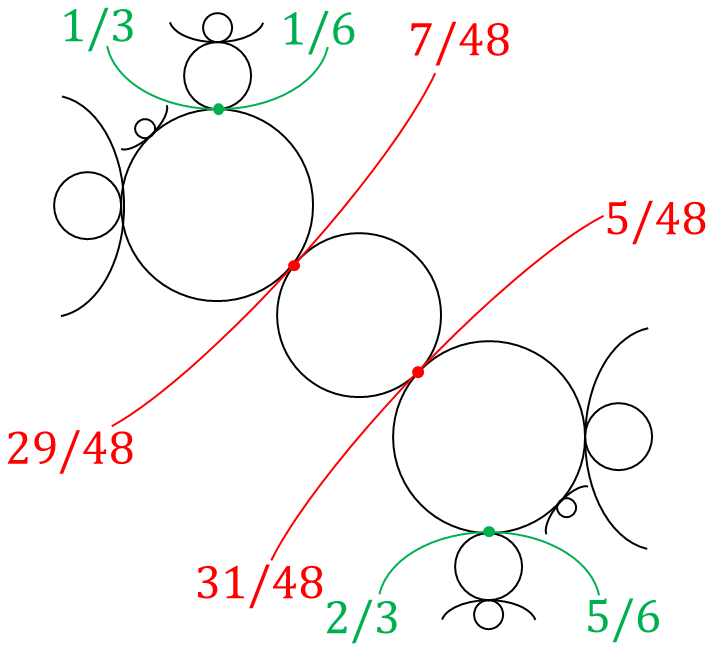}
    \end{minipage}
    \begin{minipage}{0.28\textwidth}
        \includegraphics[width=\textwidth]{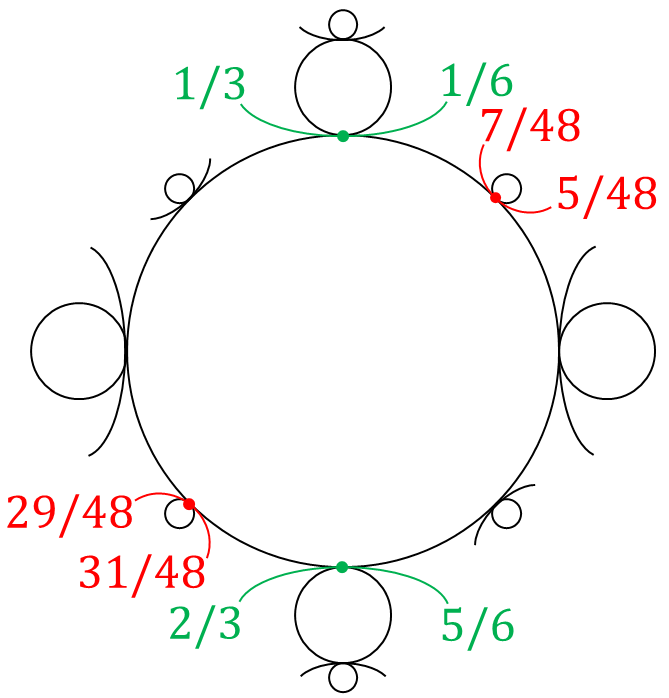}
    \end{minipage}
    \begin{minipage}{0.38\textwidth}
        \includegraphics[width=\textwidth]{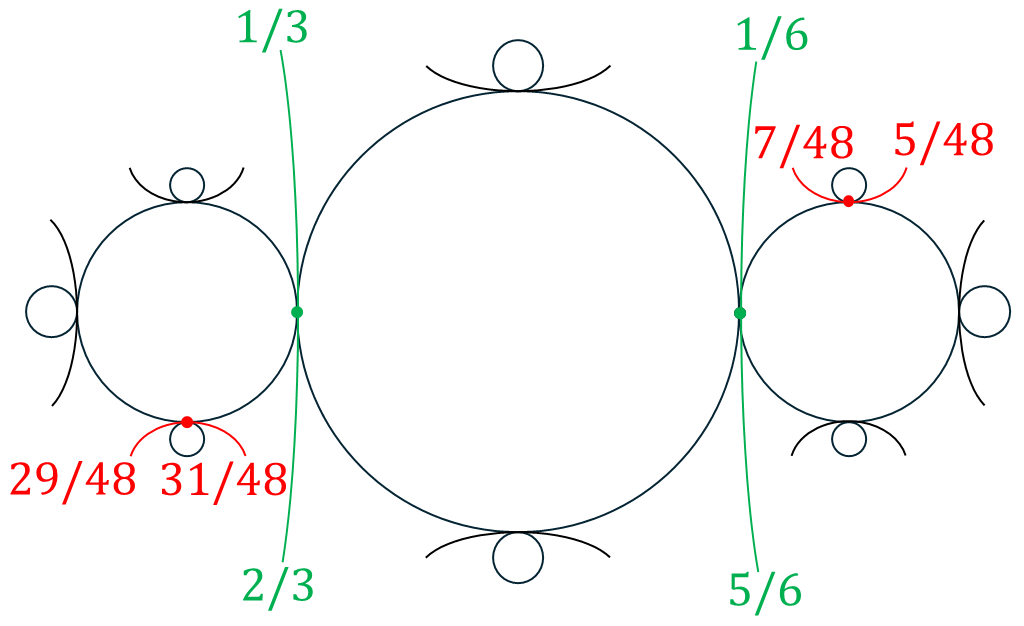}
    \end{minipage}
\caption{\label{fig:type2Jsteps}  Chain of two steps of Fatou component changes for a \textbf{Type $N=2$} example, when $v_-\in T$. Leftmost represents $\Jnot$, center is the intermediate stage and right is $\Jminus$. }
\end{figure}

\begin{figure}
        \includegraphics[width=0.6\textwidth]{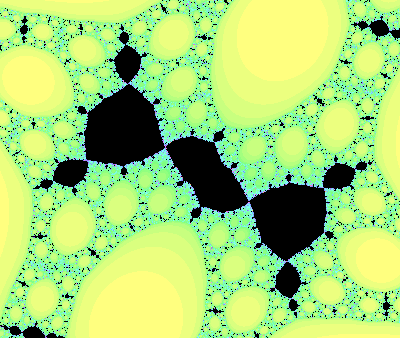}
\caption{\label{fig:Type2} For $\Rnab$ with $n=3$, $a=0.0539-0.0118i$, $b=0.01+0.03i,$ it appears $v_-$ lies in a component 
$T=
\FClabel{5/24}{11/48}{13/48}{7/24}$ in $\Jminus$, resulting in an an altered \textbf{Type 2} $J_0$. 
} 
\end{figure}

\begin{example} (\textbf{Type 2})
\label{example:Type2}
For an example of type $N=2$, consider the case where $v_- \in \FC_-$ with the identified angle name $T=\FClabel{{5}/{24}}{{11}/{48}}{{13}/{48}}{{7}/{24}}$. We illustrate this example in Figures~\ref{fig:type2lamsteps} and \ref{fig:type2Jsteps}.
The preimages of $\FC_c$ under $P_c$ are $RT=\FClabel{5/48}{11/96}{13/96}{7/48}$ (``right top'') and $LB=\FClabel{29/48}{59/96}{61/96}{31/48}$ (``left bottom''). Since $T$ in $\Kminus$ contains a critical value, $RT$ and $LB$ need to be combined into one preimage component in $\Jnot$ which contains a critical point, but at this stage $RT$ and $LB$ share no neighboring component. Hence they cannot split and reidentify their angles because there is no single component that can be split apart allowing $RT$ and $LB$ to combine. 

However, if we imagine moving $v_-$ from $L$ into $M=
\FClabel{1/6}{1/3}{2/3}{5/6}$ on its way toward $T$, taking each component traveled as a distinct step, we see that in fact we need two pairs of angle identifications to split and reidentify. When $v_-$ passes into $M$, we change the identifications from $\frac{1}{6} \sim \frac{5}{6}$, $\frac{1}{3} \sim \frac{2}{3}$ to $\frac{1}{6} \sim \frac{1}{3}$, $\frac{2}{3} \sim \frac{5}{6}$ as detailed in Type 1-1, Proposition~\ref{prop:type1M}. Recall that this reidentification causes the two preimages $L$ and $R$ of $M$ to combine into one component. 
Specifically, $L$ and $R$ have combined to form a new component 
$\FClabel{1/12}{5/12}{7/12}{11/12}$, which now lies in the center of $\Jnot$.
Then $RT$ and $LB$ are now each adjacent to this new component, %$\FC_{LR}$, 
and when we split and reidentify angles to achieve 
$\frac{5}{48} \sim \frac{31}{48}$ and $\frac{7}{48} \sim \frac{29}{48}$, it is this just formed component that splits back apart so that $\Jnot$ still maps $2:1$ onto $\Jminus$; that is, $\FClabel{1/12}{5/12}{7/12}{11/12}$ is split by the reidentifications $\frac{5}{48} \sim \frac{31}{48}$ and $\frac{7}{48} \sim \frac{29}{48}$ into two components $
\FClabel{1/12}{5/48}{31/48}{11/12}$ and $\FClabel{7/48}{5/12}{7/12}{29/48}$, and $RT$ and $LB$ combine to make one component 
$\FClabel{5/48}{7/48}{29/48}{31/48}$. 

Figure~\ref{fig:type2lamsteps} shows the chain of lamination changes, and Figure~\ref{fig:type2Jsteps} the chain of Fatou component changes. 
Figure~\ref{fig:Type2} is an $\Rnab$ which appears to be \textbf{Type 2}. 
\end{example} 

Now we are ready to state and prove the general case of $v_-$ being $N$ Fatou components away from $L$. 

\begin{theorem} \label{thm:typeN} (\textbf{Type $N \geq 1$})
    Suppose that $v_-$ lies in a Fatou component $\FC_-$ of $\Kminus$ which shares an identified angle name with the Fatou component $\FC_c$ of $K_c$ such that $\FC_c$ is any Fatou component other than 
    $L=\FClabel{1/3}{5/12}{7/12}{2/3}$. 

    Starting with $\FC_0=L$ and ending with $\FC_N=\FC_c$, let $\FC_0, \dots, \FC_N$ be the shortest possible path of adjacent Fatou components in $K_c$ from $L$ to $\FC_c$. For $i=1, \dots, N$, let $\FC_i$'s preimage Fatou components be labeled $\FClabel{a_1^i}{b_1^i}{a_2^i}{b_2^i}$ and $\FClabel{a_3^i}{b_3^i}{a_4^i}{b_4^i}$ , noting that $a_1^i \sim b_2^i, b_1^i \sim a_2^i, a_3^i \sim b_4^i, b_3^i \sim a_4^i$. 

    Then the lamination $\lam_0$ for $\gamma_0$ is the same as the lamination $\lam_-$ for $\gamma_-$
    except that for each $i=1, \dots, N$, the identifications change to $a_1^i \sim b_4^i$ and $b_2^i \sim a_3^i$ (with the identifications $b_1^i \sim a_2^i, b_3^i \sim a_4^i$ unchanged), with one exception: if $\FC_1 = M=
    \FClabel{1/6}{1/3}{2/3}{5/6}$, 
    then $\lam_0$ has the identified angles $b_1^1=\frac{1}{6} \sim \frac{1}{3}=a_3^1$ and $b_4^1=\frac{2}{3} \sim \frac{5}{6}=a_2^1$.
    
    Furthermore, these angle identification changes result in the lamination diagram $\lam_0$ for $\gamma_0$ differing from the lamination diagram $\lam_-$ for $\gamma_-$ by $2N+1$ leaves.
\end{theorem}

\begin{proof}

First, we provide an explanation of why an iterative process is needed.

We note that we let $\FC_0, \dots, \FC_N$ be the shortest path from $L$ to $\FC_c$ because any longer path would involve backtracking. Then any step taken off the most direct path would be undone, and therefore would have no affect on the final lamination. Hence without loss of generality, we consider the shortest possible path.

Recall by Corollary~\ref{cor:lamination_symmetry}, lamination diagrams for all $\gamma_{m,j}$'s have $180^{\circ}$ rotational symmetry. 
So, there are two chords, opposite each other in the lamination, whose pairs of identified angles represent two points in $\CC$ on the boundaries of two preimage Fatou components of $\FC_c$.

Observe that the \textbf{central leaf} of the lamination diagram for $J_-$ corresponds to $M$ and is the one which contains $0 \in \mathbb{D}$. 
As established in Propositions~\ref{prop:type1M} and \ref{prop:type1notM}, to combine two Fatou components into one which contains the critical point (which corresponds to $0$ in the lamination disk), hence a ``new'' central leaf, we need the two chords that bound the corresponding leaves in the lamination diagram to be broken, and then glued together in a different way, like closing and then opening a paper fortune teller the other way. However, unless those two Fatou components that will join are adjacent to $M$, we can't just make that one change in the lamination diagram, as it would create chord intersections, thus an invalid lamination. 
Thus combining the two opposite leaves requires a series of steps, changing identifications and hence leaves/Fatou components that are essentially ``in the way'' as we walk $v_-$ from its ``expected'' position in $L$ to its actual position in $\FC_-$. At each step, one pair of angle identifications is split and re-identified, similar to what occurs for \textbf{Type 1} described above. 
We present this as a finite induction. 

\textbf{Base Step:} If $\FC_c$ is any Fatou component of $K_c$ other than $L$, then $\FC_1$ is a component adjacent to $L$. Then $\FC_1$ is either $M$ and we create an (intermediate if $N>1$) lamination $\lam_{N-1}$ via the changes detailed in Proposition~\ref{prop:type1M}, or $\FC_1$ is any component adjacent to $L$ other than $M$ and the changes detailed in Proposition~\ref{prop:type1notM} must happen. In either case, $\lam_{N-1}$ differs from $\lam_-$ by $3 = 2+1$  leaves. 
If $N=1$ then $\lam_{N-1}=\lam_0$, and the differences between $\lam_0$ and $\lam_-$ as listed in the theorem statement have occurred, consistent with Propositions~\ref{prop:type1M} and~\ref{prop:type1notM}. 

\textbf{Inductive Step:} Suppose now that we have split and re-identified a pair of angle chords $k$ times, for some $k\in\{1,\ldots,N-1\}$,
as we imagined moving $v_-$ from $L=\FC_0$ through $k$ adjacent Fatou components, $\FC_1, \ldots, \FC_{k}$, and now we have an intermediate lamination diagram $\lam_{N-k}$, which differs from $\lam_c$ by $2k+1$ leaves.  
Now our imagined position of $v_-$ is in $\FC_{k}$. 
We let the two preimage Fatou components of $\FC_{k+1}$ be labelled $\FClabel{a_1^{k+1}}{b_1^{k+1}}{a_2^{k+1}}{b_2^{k+1}}$ and $\FClabel{a_3^{k+1}}{b_3^{k+1}}{a_4^{k+1}}{b_4^{k+1}}$. Based on the change of components carried out in steps $i=1, \dots, k$, we have that $a_1^{k+1} \sim b_2^{k+1}$ and $a_3^{k+1} \sim b_4^{k+1}$ are the places where these Fatou components meet the current central component of the intermediate $\lam_{N-k}$, which we'll label as $\FClabel{a}{b}{c}{d}$. Recall again that $a_1^i$ and $a_3^i$ double to the same angle, and so do $b_2^i$ and $b_4^i$. Then to combine these two components into one new central component for our next intermediate step, but respecting the dynamics of the map and the fact that the angle assignments $\gamma_0$ respect the dynamics as described in Theorem~\ref{thm:angle assignments}, 
we break the identifications so that  $a_1^{k+1} \nsim b_2^{k+1}$ and $a_3^{k+1} \nsim b_4^{k+1}$ and we create the new identifications $a_1^{k+1} \sim b_4^{k+1}$ and $b_2^{k+1} \sim a_3^{k+1}$, maintaining all other identifications. Then the new central Fatou component which must contain the critical point is $\FClabel{a_1^{k+1}}{b_2^{k+1}}{a_3^{k+1}}{b_4^{k+1}}$ and the previous central Fatou component has been split apart into $\FClabel{a}{a_1^{k+1}}{b4^{k+1}}{d}$ and $\FClabel{b_2^{k+1}}{b}{c}{a_3^{k+1}}$. This creates the next intermediate lamination $\lam_{N-(k+1)}$. 

Note since the ``new'' central component from the prior step has been changed again, that change doesn't increase the number of leaves in $\lam_{N-(k+1)}$ different from $\lam_{N-k}$, so the number of leaves in $\lam_{N-(k+1)}$ different from $\lam_0$ is $2k+1+2 = 2(k+1)+1.$

The induction is finite: if $k+1=N$, the lamination $\lam_{N-N}$ that was just created is not intermediate but is $\lam_0,$ and the total number of different leaves is $2N+1$. 
\end{proof}

\begin{figure}[h]
    \begin{minipage}{0.4\textwidth}
        \includegraphics[width=\textwidth]{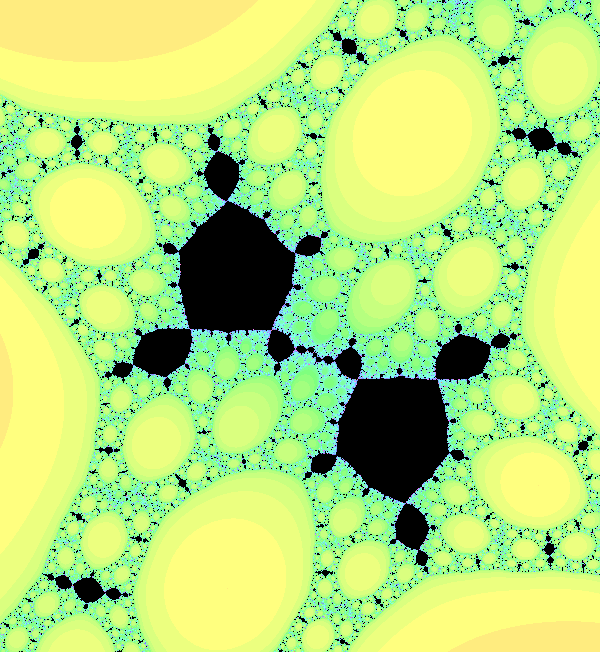}
    \end{minipage}
    \begin{minipage}{0.45\textwidth}
        \includegraphics[width=\textwidth]{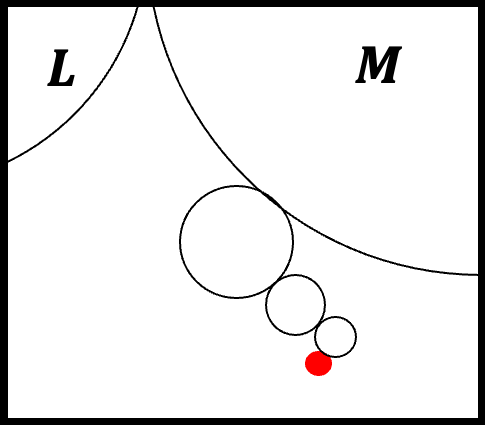}
    \end{minipage}
    \caption{\label{fig:Type4} A $\Jnot$ example of \textbf{Type $N=5$} is shown on the left, using $n=3$, $a \approx 0.054297-0.012066i$, $b=0.01+0.03i$. On the right is a representation of the placement of $v_-$ within $\Jminus$, $v_-$ lies in the component shaded in red. These Fatou components are not to scale, as each actual Fatou component along this path is much smaller than the prior one, so we enlarged the smaller components in order to be clearly visible.} 
\end{figure}

\begin{example}  \label{example:Type5}
See Figure~\ref{fig:Type4} for what appears to be a \textbf{Type 5} example, with a computer generated image of $\Jnot$ and a diagram indicating the location of $v_-$. 
\end{example}

\smallskip

\noindent \textbf{Observation on Fatou component size.} We observe that the alterations of each $\Jnot$ did not only consist of changing angle assignments, but there's also a pattern regarding the relative sizes of components. Note in typical polynomial Julia sets, the component of $K$ containing the critical point is large relative to the size of the critical value component (in other words, the sole preimage component of the critical value component is ``large'', relatively). Since the derivative at the critical point is zero, naturally if a component in a baby $J$ preimage  happens to contain the other critical value $v_-$, its preimage component, which lies in an ``altered'' baby $J$, contains a critical point, and as such is larger than would have been expected.
However, unless $\FC_-=M$, new Fatou components are formed from merging together what were smaller Fatou components, so the visually central components may become smaller.

%-----------------------------------------
\subsection*{Future Explorations} 
\begin{figure}
    \includegraphics[width=0.6\textwidth]{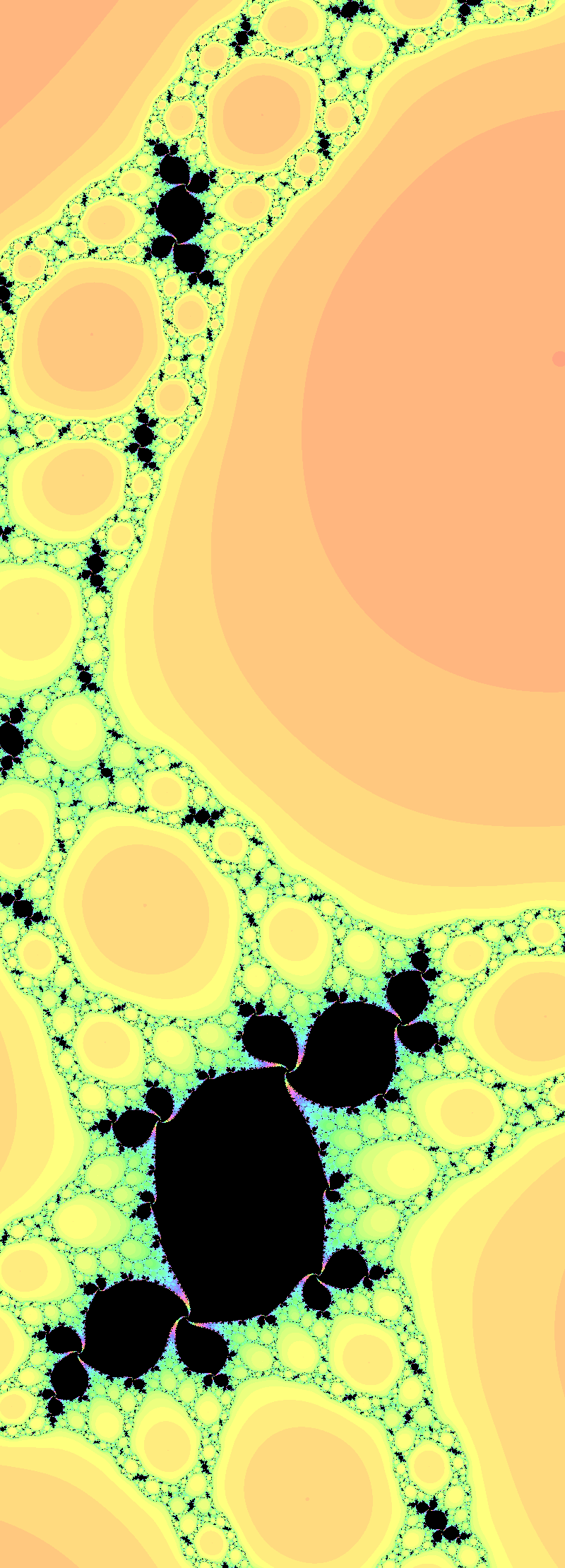}
    \caption{\label{fig:altered rabbits} Observe both typical looking and ``altered'' rabbits - note the mismatched ``ear'' size on the larger ``altered'' rabbit. Here $n=3$, $a\approx 0.106497+0.077201i$, $b=0.01+0.05i$.}
\end{figure}

We restricted to the basilica case for ease of exposition and illustration. However, this phenomena is easy to observe for other quadratic Julia sets (Douady rabbits, aeroplanes, etc.), and a more general result could be interesting future work. 
Figure~\ref{fig:altered rabbits} shows portions of a Julia set of an $\Rnab$ which appears to have altered rabbit preimages of a standard baby rabbit-like Julia set. 

Finally, in previous work, Boyd and co-authors \cite{boyd_mitchell, BoydHoeppner1} applied Douady and Hubbard's results to prove the existence of baby Julia sets for $P_c$ in the Julia set for several different restrictions on the parameters $n,a,b$ of $\Rnab$, and baby Mandelbrot sets in parameter spaces.  So we know that Assumptions A3 and A4 actually occur. In \cite{boyd_mitchell}, Section~5.3, it is shown that for $n\geq 11$, as $a$ varies from $0.1$ to $1$, there is an area of $b$-values for which both critical orbits are bounded and each is creating a baby quadratic-like portion of the Julia set of $\Rnab$ associated with it. Nearby this area is where we can experimentally observe many examples of Assumptions A3 and A5 appearing together (for example, $\Rnab$ for $n=11, a=0.2075+0.0i,$ and $b=-0.004483852355144613+0.0032857624406795257i$ appears to be of \textbf{Type 1-1} and lies in a region in which we proved Assumption A3 holds), helping to guide our explorations of what the dynamics must be if this does occur. However, actually establishing that there are parameters for which $\Rnab$ satisfies Assumptions A1 through A5 would be a good avenue for future study.

%=========== BIBLIOGRAPHY==============
%\bibliographystyle{alpha} 

\end{document}